\documentclass{amsart}

\usepackage{amsmath}
\usepackage{amssymb}
\usepackage{datetime}
\usepackage{enumitem}
\usepackage{hyperref}

\newtheorem{theorem}{Theorem}[section]
\newtheorem{lemma}[theorem]{Lemma}
\newtheorem{proposition}[theorem]{Proposition}
\newtheorem{corollary}[theorem]{Corollary}
\newtheorem*{main}{Theorem A}
\newtheorem*{main2}{Theorem B}
\newtheorem*{main3}{Theorem C}

\theoremstyle{definition}

\theoremstyle{remark}
\newtheorem{remark}[theorem]{Remark}

\numberwithin{equation}{section}

\newcommand{\R}{\mathbb{R}}
\newcommand{\C}{\mathbb{C}}

\newcommand{\Kill}{\mathrm{Kill}}
\newcommand{\Iso}{\mathrm{Iso}}

\newcommand{\OO}{\mathcal{O}}
\newcommand{\HH}{\mathcal{H}}
\newcommand{\GG}{\mathcal{G}}
\newcommand{\VV}{\mathcal{V}}
\newcommand{\WW}{\mathcal{W}}

\newcommand{\FF}{\mathcal{F}}

\newcommand{\g}{\mathfrak{g}}
\newcommand{\h}{\mathfrak{h}}
\newcommand{\s}{\mathfrak{s}}

\newcommand{\so}{\mathfrak{so}}
\newcommand{\su}{\mathfrak{su}}

\newcommand{\sli}{\mathfrak{sl}}
\newcommand{\gl}{\mathfrak{gl}}

\newcommand{\rad}{\mathrm{rad}}

\newcommand{\LA}{\langle}
\newcommand{\RA}{\rangle}
\newcommand{\TG}{\widetilde{G}}
\newcommand{\TM}{\widetilde{M}}
\newcommand{\TN}{\widetilde{N}}
\newcommand{\TSL}{\widetilde{\SL}}
\newcommand{\sS}{\mathrm{S}}
\newcommand{\SL}{\mathrm{SL}}
\newcommand{\Rank}{\text{rank}}

\newcommand{\W}{\mathrm{W}}

\begin{document}

% Title Page
\title{Rigidity of an Isometric $\SL(3,\R)$-Action}
\author{Ra\'ul Quiroga-Barranco}
\address{Centro de Investigaci\'on en Matem\'aticas, Apartado Postal 402, Guanajuato, Guanajuato, 36250, Mexico}
\email{quiroga@cimat.mx}
\author{Eli Roblero-M\'endez}
\address{Department of Mathematics, 346 Lockett hall, Louisiana State University, Baton Rouge, LA, 70803, USA}
\email{elir@math.lsu.edu}

%\date{March 04, 2016}

\begin{abstract}
 We characterize the universal covering of connected analytic pseu\-do-Riemannian manifolds which admit a
 non-trivial and isometric action of the simple Lie group $\SL(3,\R)$ with a dense orbit preserving a
 finite volume. If such manifold is also weakly irreducible we prove that $M$ is isometric to, or a quotient
 space of, a simple Lie group containing $\SL(3,\R)$.
\end{abstract}

\maketitle

\section*{Introduction}

 Let $G$ be a connected non-compact simple Lie group acting on a connected, analytic manifold $M$ preserving
 a pseudo-Riemannian structure. It is conjectured that such actions are rigid in the sense that they
 restrict the possibilities of the manifold. It is expected that any such actions, under non-trivial
 conditions, must be an algebraic double coset. That is, $M\cong K\backslash H/\Gamma$, such that $H$ is
 a Lie group and there is a group homomorphism $G\hookrightarrow H$ whose image commutes with $K$, a compact
 subgroup of $H$, and $\Gamma\subset H$ is a lattice. Therefore, the $G$-action is given by a natural left
 action on $M\cong K\backslash H/\Gamma$. Some results in this direction are obtained in \cite{OQ} and
 \cite{OQ-Upq}, proving that with some extra geometric conditions such $G$-actions imply that the manifolds
 are of double-coset type.

 In \cite{OQ}, for $M$ a complete and weakly irreducible manifold with a non-transitive $G$-action and with
 a dense orbit, it is proved that the dimension of $M$ has a lower bound, in terms of the theoretical
 properties of the corresponding representation of $\g$, the Lie algebra of $G$. In particular, \cite{OQ}
 explains in detail the case $G=\widetilde{SO}_0(p,q)$, with $p+q\geq4$, determining the manifold $M$.

 The paper applies the techniques from \cite{OQ} to the case $G=\SL(3,\R)$, a connected,
 simple Lie group with Lie algebra $\sli(3,\R)$. In our case, $\SL(3,\R)$ acts isometrically on $M$, a
 complete, pseudo-Riemannian manifold such that $8<\dim(M)\leq14$ with a dense orbit. The assumption on the
 dimension of $M$ eliminates the condition of non-transitivity of the action.

 In \cite{OQ} and \cite{OQ-Upq} we observe the study of irreducible representations of groups
 preserving a non-degenerate, symmetric, bilinear form. Such study is an important tool to understand the
 normal bundle of the foliation generated by the action of the group. That understanding together with
 the properties of the tangent bundle of the foliation give us information which restrict the possibilities of
 the manifold $M$.

 Here, we analyze the representation of $\SL(3,\R)$ of minimal dimension with the property of preserving a
 non-degenerate, symmetric bilinear form, which is non-irreducible. Hence, this paper generalizes such
 corresponding part of the previous works.

 As in \cite{OQ}, we obtain a lower bound on the dimension of $M$ which strongly determines
 $\TM$, the universal covering of $M$. This is the result of our main theorem.

 \begin{main}\label{th: main theorem}
  Let $M$ be a connected analytic pseudo-Riemannian manifold. Suppose that $M$ is complete, has finite volume
  and admits an analytic and isometric $SL(3,\R)$-action with a dense orbit. If $8<\dim(M)\leq14$, then $\TM$
  is isomorphic to one, and only one, of the following:
  \begin{enumerate}
   \item[$(i)$] $\TSL(3,\R)\times\TN$, where $\TN$ is a complete pseudo-Riemannian manifold.
   \item[$(ii)$] $G_{2(2)}$, the simply connected real Lie group related to the non-compact real form of
    the exceptional simple Lie algebra $\g_2^\C$.
   \item[$(iii)$] $\R\backslash\TSL(4,\R)$.
  \end{enumerate}
 \end{main}

 Recall that a connected pseudo-Riemannian manifold is \emph{weakly irreducible} if there is no proper,
 non-degenerate invariant subspace of the tangent space at some (at hence any) point invariant under the
 restricted holonomy group at that point. If in the hypothesis of Theorem A we assume that the manifold
 $M$ is weakly irreducible then case (i) is no allowed, that is:

 \begin{main2}\label{th: Main Theorem2}
  Let $M$ be a connected analytic pseudo-Riemannian manifold. Suppose that $M$ is complete, has finite volume
  and admits an analytic and isometric $SL(3,\R)$-action with a dense orbit. If $8<\dim(M)\leq14$ and $M$ is a
  weakly irreducible manifold, then $\TM$ is isomorphic to one of the following:
  \begin{enumerate}
   \item[$(1)$] $G_{2(2)}$.
   \item[$(2)$] $\R\backslash\TSL(4,\R)$.
  \end{enumerate}
  The corresponding isomorphism is $\TSL(3,\R)$-equivariant, where the $\TSL(3,\R)$-action on it is
  induced by some non-trivial homomorphism of $\TSL(3,\R)$ into $G_{2(2)}$ or $\TSL(4,\R)$, respectively.
  We can also rescale the metric along
  $\TSL(3,\R)$-orbits  and their normal bundle to assume that such isomorphism is a local isometry for the
  bi-invariant pseudo-Riemannian metric on $G_{2(2)}$ or on $\R\backslash\TSL(4,\R)$ is given by the Killing
  form of $\g_{2(2)}$ or $\sli(4,\R)$, respectively.
 \end{main2}

 Aditionally, if we assume that the universal covering of the manifold is not isometric to a quotient space
 then the option (1) can be eliminated in the previous theorem. This, $\TM\cong G_{2(2)}$ and $M$ is isometric
 to a quotient of $G_{2(2)}$ over a discrete subgroup. Which is the result of the next theorem.

 \begin{main3}\label{th: Main Theorem3}
  Let $M$ be a connected analytic pseudo-Riemannian manifold. Suppose that $M$ is complete, has finite volume
  and admits an analytic and isometric $SL(3,\R)$-action with a dense orbit. If $8<\dim(M)\leq14$ and $M$ is a
  weakly irreducible manifold such that $\TM$ is not isomorphic to a quotient map then there exist
  \begin{itemize}
   \item a lattice $\Gamma\subset G_{2(2)}$, and
   \item an analytic finite covering map $\varphi:G_{2(2)}/\Gamma\to M$
  \end{itemize}
  such that $\varphi$ is $\TSL(3,\R)$-equivariant map, where the $\TSL(3,\R)$-action on $G_{2(2)}/\Gamma$ is
  induced by some non-trivial homomorphism $\TSL(3,\R)\to G_{2(2)}$. We can also rescale the metric along
  $\TSL(3,\R)$-orbits  and their normal bundle to assume that $\varphi$ is a local isometry for the
  bi-invariant pseudo-Riemannian metric on $G_{2(2)}$ given by the Killing form of its Lie algebra.
 \end{main3}

 The proofs of Theorems A, B and C are based on the application of representation theory to the Killing
 vector fields that centralize the action of the group $\TSL(3,\R)$. A fundamental result for this work
 is Proposition \ref{prop: g(x)}, which appears in \cite{OQ} and \cite{Q-Ann} and whose detailed proof
 can be found in \cite{Q-TINC}. Such proposition shows the existence of a subalgebra of Killing vector fields
 on $\TM$ which vanish at some fixed point. The previous algebra is used to show properties of the centralizer
 of the $\TSL(3,\R)$-action, that we denote by $\HH$. The analysis of such properties of $\HH$ lead us to
 identify the possible options of $\TM$.

 The organization of this paper is the following: In Section \ref{sect: Action} we present the minimal
 dimension representation $W$ of $\sli(3,\R)$ preserving a metric and therefore the decomposition of $\so(W)$
 as a direct sum of irreducible $\sli(3,\R)$-modules. In Section \ref{sect: Centralizer}, we show results
 that guarantee the existence of the centralizer of the action on the manifold $M$. In Section
 \ref{sect: Centralizer II} we study the properties of the centralizer as a $\sli(3,\R)$-module which
 lead us to obtain a lower bound of the dimension of $M$. In Section \ref{sect: Structure}, we analyze
 the structure of the centralizer that permits to restrict the possibilities of the manifold $M$.
 Finally, in Section \ref{sect: Proof}, we use the previous result to prove Theorems A, B and C.

\section{Autodual non-trivial Representations of \texorpdfstring{$\sli(3,\R)$}{}}\label{sect: Action}

 Recall that $\SL(3,\R)$ denotes the \textbf{special linear group} of degree $3$ over the field of real numbers.
 The group $\SL(3,\R)$ is a connected, non-compact, real simple Lie group and Lie algebra denoted as
 $\sli(3,\R)$.

 Let $\rho$ be a representation of $\sli(3,\R)$ on the real vector space $V_0$ which
 induces a complex representation $\rho(\C):\sli(3,\C)\to V$, where $V=V_0\otimes\C$. If a non-trivial
 representation $\rho$ preserves a non-degenerated symmetric bilinear form then  the representation
 $\rho(\C)$ is \emph{autodual}.

 By the Weyl Character Formula, the lower dimensions of non-trivial irreducible representations of
 $\sli(3,\C)$ are $3$, $6$ and $8$. Where $3$ is the dimension of the natural representation of
 $\SL(3,\C)$ onto $\C^3$ or onto its dual ${\C^3}^*$, which are non-autodual. Therefore, one can see
 that the non-trivial representations of $\sli(3,\R)$ with dimension less than $6$ are non-autodual.

 The previous result is useful to show that the minimum dimension of a non-trivial representation of
 $\sli(3,\R)$ is $6$, which is a corollary of the next lemma.

 \begin{lemma}\label{lem: Bil Form}
  The minimal non-trivial real representation of $\sli(3,\R)$ preserving a non-degenerated bilinear symmetric form is
  isomorphic to $\R^3\oplus\R^{3*}$. That bilinear form has signature $(3,3)$.
 \end{lemma}
 \begin{proof}
  Let $\rho:\sli(3,\R)\to V_0$ be a non-trivial real representation preserving a non-degenerated symmetric bilinear
  form $\LA\cdot,\cdot\RA_0$. By the above we have that $\dim_\R(V_0)\geq6$.

  First, we assume that $\dim_\R(V_0)=6$. If $\rho$ is an irreducible representation then $\rho(\C)$ is a
  complex irreducible representation. Therefore, Section 13 of \cite{Fulton} implies
  $V\simeq\text{Sym}^2(\C^3)$ or $V\simeq\text{Sym}^2(\C^{3*})$. In both cases, the representation is non-autodual.

  On the other hand, if $\rho$ is reducible (and therefore $\rho(\C)$ is also reducible) then $V$ is isomorphic
  to one of the following vector spaces: $\bigoplus_{j=1}^6\C$, $\C^3\oplus\bigoplus_{j=1}^3\C$,
  $\C^{3*}\oplus\bigoplus_{j=1}^3\C$, $\C^3\oplus\C^3$, $\C^{3*}\oplus\C^{3*}$ or $\C^3\oplus\C^{3*}$.

  By the properties of representation of $\sli(3,\C)$ on
  $\C$, $\C^3$ and $\C^{3*}$ we have that $V$ cannot be isomorphic, as a $\sli(3,\C)$-module, to
  $\bigoplus_{j=1}^6\C$, $\C^3\oplus\bigoplus_{j=1}^3\C$, $\C^{3*}\oplus\bigoplus_{j=1}^3\C$, $\C^3\oplus\C^3$ and
  $\C^{3*}\oplus\C^{3*}$. Hence, $V_0$ cannot be isomorphic (as $\sli(3,\R)$-module) to
  $\bigoplus_{j=1}^6\R$, $\R^3\oplus\bigoplus_{j=1}^3\R$, $\R^{3*}\oplus\bigoplus_{j=1}^3\R$, $\R^3\oplus\R^3$ and
  $\R^{3*}\oplus\R^{3*}$.

  Next, we prove the existence of a non-degenerated symmetric bilinear form $\LA\cdot,\cdot\RA_0$ in
  $\R^3\oplus\R^{3*}$ which is $\sli(3,\R)$-invariant. Let $v,v'\in\R^3\oplus\R^{3*}$ be given, then there exist unique elements $p,p'\in\R^3$,
  $q,q'\in\R^{3*}$ such that $v=(p,q)$ and $v'=(p',q')$. We define $\LA v,v'\RA_0$ as follows
  \begin{equation*}
   \LA v,v'\RA_0=\LA(p,q),(p',q')\RA_0:=q(p')+q'(p),
  \end{equation*}
  where $q(p)$ is the evaluation of the element $q$ in the vector $p$. Note that
  $\LA v,v'\RA_0$ is a non-degenerated symmetric non-degenerated bilinear form. Let $A$ be an
  arbitrary but fixed element of $\sli(3,\R)$ then,
  \begin{eqnarray*}
    \LA A\cdot v,v'\RA+\LA v,A\cdot v'\RA_0 &=& \LA A\cdot(p,q),(p',q')\RA_0+\LA(p,q),A\cdot(p',q')\RA_0 \\
     &=& \LA(A\cdot p,A\cdot q),(p',q')\RA_0+\LA(p,q),(A\cdot p',A\cdot q')\RA_0 \\
     &=& (A\cdot q)(p')+q'(A\cdot p)+q(A\cdot p')+(A\cdot q')(p)\\
     &=& q(-A\cdot p')+q'(A\cdot p)+q(A\cdot p')+q'(-A\cdot p)\\
     &=& 0,
  \end{eqnarray*}
  which proves that $\LA\cdot,\cdot\RA_0$ is $\sli(3,\R)$-invariant.

  Now, let $\LA\cdot,\cdot\RA$ be a non-degenerated symmetric bilinear form on $\R^3\oplus\R^{3*}$ which is
  $\sli(3,\R)$-invariant. If $\{e_1,e_2,e_3\}$ denotes the canonical basis of $\R^3$ then, given
  $i,j\in\{1,2,3\}$ such that $i\neq j$, we can find an element $A_{ij}\in\sli(3,\R)$ such that
  \begin{equation*}
   A_{ij}(e_i)=e_i \qquad \text{and} \qquad A_{ij}(e_j).% \qquad \text{for}
  \end{equation*}
  The $\sli(3,\R)$-invariance of $\LA\cdot,\cdot\RA$ implies that
  \begin{equation*}
   0=\LA A_{ij}(e_i),e_j\RA+\LA e_i,A_{ij}(e_j)\RA=2\LA e_i,e_j\RA,
  \end{equation*}
  hence $\LA e_i,e_j\RA=0$. Thus, the signature of the bilinear form $\LA\cdot,\cdot\RA_0$ on
  $\R^3\oplus\R^{3*}$ is $(3,3)$.
 \end{proof}

 Since the Lie algebra $\sli(3,\R)$ preserves a non-degenerate symmetric bilinear form on $\R^3\oplus\R^{3*}$ with
 signature $(3,3)$, there exists a non-trivial Lie algebra homomorphism $\sli(3,\R)\to\so(3,3)$. The
 simplicity of $\sli(3,\R)$ implies that such homomorphism is injective and therefore $\so(3,3)$ has a structure of
 $\sli(3,\R)$-module.

 Next, we analyze the decomposition of $\so(3,3)$ into a direct sum of irreducible $\sli(3,\R)$-submodules.

 From Table II in \cite{Berger},  $\big(\so(3,3),\gl(3,\R)\big)$ is a symmetric pair where
 $\gl(3,\R)=\sli(3,\R)\oplus\R$ is its decomposition as a direct sum of irreducible $\sli(3,\R)$-modules.

 \begin{lemma}\label{lem: D irred}
  $\so(3,3)$ is isomorphic, as a $\sli(3,\R)$-module, to $\sli(3,\R)\oplus\R^3\oplus\R^{3*}\oplus\R$.
 \end{lemma}
 \begin{proof}
 In general, Table II in \cite{Berger}, shows $\big(\so(n,n),\sli(n,\R)\oplus\R\big)$ is a symmetric
 pair  such that
 \begin{equation*}
   \so(n,n)=\sli(n,\R)\oplus\R\oplus\nu_2(\sli(3,\R))
 \end{equation*}
 where $\nu_2(\sli(3,\R))$ is a self-adjoint $\sli(n,\R)$-module containing $\pi_2(\sli(n,\R))$,
 the irreducible representation of $\sli(n,\R)$ correspondent to its second fundamental weight $\varpi_2$.
 Using a similar analysis, about self-adjoint representations, of Appendix A in \cite{OQ-Upq} we have that
 $\nu_2(\sli(3,\R))=\pi_2(\sli(n,\R))\oplus\pi_{n-2}(\sli(n,\R))$.

 In the case $n=3$ we have that $\pi_2(\sli(3,\R))={\R^3}^*$ and $\pi_1(\sli(3,\R))=\R^3$.
 Therefore
 \begin{equation*}
   \so(3,3)=\sli(3,\R)\oplus\R\oplus\R^3\oplus{\R^3}^*.
 \end{equation*}
 \end{proof}

 The following corollary is a consequence of the previous lemma.
 \begin{corollary}\label{cor: subs}
  The subalgebras of $\so(3,3)$ that are at the same time $\sli(3,\R)$-submodules, with the structure
  of module induced by an
  injection of $\sli(3,\R)$ into $\so(3,3)$, are isomorphic to one of the following: $0$, $\sli(3,\R)$, $\R^3$,
  $\R^{3*}$, $\R$, $\sli(3,\R)\oplus\R^3$, $\sli(3,\R)\oplus\R^{3*}$, $\sli(3,\R)\oplus\R$, $\R^3\oplus\R$,
  $\R^{3*}\oplus\R$, $\sli(3,\R)\oplus\R^3\oplus\R$, $\sli(3,\R)\oplus\R^{3*}\oplus\R$ or to $\so(3,3)$.
 \end{corollary}
 \begin{proof}
   We only recall that $\big(\so(3,3),\sli(3,\R)\oplus\R\big)$ is a symmetric pair.
 \end{proof}

\section{Isometric actions of the simple Lie group \texorpdfstring{$\SL(3,\R)$}{}}\label{sect: Centralizer}

 In this section we assume $G$ is a connected non-compact simple Lie group with Lie algebra $\g$ and $M$ a connected
 and analytic finite-volume pseudo-Riemannian manifold where $G$ acts analytic and isometrically with a dense orbit.

 Since every isometric $G$-action on a manifold $M$ with a dense orbit is locally free (see \cite{Szaro}),
 the orbits of the $G$-action define a foliation that we denote by $\FF$. Then, for every $x\in M$ there exist a
 vector subspace of $T_xM$ such that $T_xM=T_x\FF\oplus T_x\FF^\perp$. The collection of these vector subspaces
 form a distribution on $M$ that we will denote as $T\FF^\perp$.

 On the other hand, the tangent bundle to the foliation $\FF$ is a trivial vector bundle isomorphic to
 $M\times\g$ under the isomorphism $M\times\g\to T\FF$, given by $(x,X)\mapsto X^*_x$. That also defines an
 isomorphism fiber of $T_x\FF$ with $\g$. For the rest of the paper, for an element $X$ in the Lie algebra of
 a group acting on a manifold, we denote by $X^*$ the vector field on the manifold whose one-parameter group of
 diffeomorphism is given by $(exp(tX))_t$ through the action on the manifold.

 The space of Killing fields of a geometric structure $\omega$ in a manifold $M$ is denoted by $\Kill(M,\omega)$,
 and $\Kill_0(M,\omega,x)$ will denote the subspace of $\Kill(M,\omega)$ consisting of vector fields that vanish
 on $x$. Here, we denote by $\sigma$ the geometric structure of the pseudo-Riemannian metric on $M$.
 Unless otherwise is indicated, we will also omit the symbol that denotes the structure of
 pseudo-Riemannian metric, in particular, $\Kill(M):=\Kill(M,\sigma)$.

 Let $V$ be a vector space, we denote by $\so(V)$ the Lie algebra of linear maps on $V$ that are
 skew-symme\-tric with respect to a non-degenerate symmetric bilinear form on $V$.

 As an immediate use of the previous definition and as consequence of the Jacobi identity we have the next
 lemma, showed in \cite{OQ}.
 \begin{lemma}\label{lem: homom to so}
  Let $N$ be a pseudo-Riemannian manifold and $n\in N$. Then, the map $\lambda_n:\Kill_0(N,n)\to\so(T_nN)$ given by
  $\lambda_n(Z)(w)=[Z,W]_n$, where $W$ is any vector field such that $W_n=w$, is a well defined homomorphism of
  Lie algebras.
 \end{lemma}

 Here, the universal covering of any manifold $N$ it will be denoted as $\widetilde{N}$.
 Next, we present a result, proved in \cite{Q-TINC}, which is fundamental in the present work.
 \begin{proposition}[{\cite[Proposition 2.3]{Q-TINC}}]\label{prop: g(x)}
  Let $G$ be a connected non-compact simple Lie group acting isometrically and with a dense orbit on a connected
  finite volume pseudo-Riemannian manifold $M$. Consider the $\TG$-action on $\TM$, lifted from the $G$-action
  on $M$. Assume that $M$ and the $G$-action on $M$ are both analytic. Then, there exists a conull subset
  $S\subset\TM$ such that for every $x\in S$ the following properties are satisfied:
  \begin{itemize}
   \item[1.] There is a homomorphism $\rho_x:\g\to\Kill(\TM)$ which is an isomorphism onto its image
    $\rho_x(\g)=\g(x)$.
   \item[2.] $\g(x)\subset\Kill_0(\TM,x)$, i.e. every element of $\g(x)$ vanishes at x.
   \item[3.] For every $X,Y\in\g$ we have
    \begin{equation*}
     [\rho_x(X),Y^*]=[X,Y]^*=-[X^*,Y^*].
    \end{equation*}
    In particular, the elements in $\g(x)$ and their corresponding local flows preserve both $\FF$ and $T\FF^\perp$.
   \item[4.] The homomorphism of Lie algebras $\lambda_x\circ\rho_x:\g\to\so(T_x\TM)$ induces a $\g$-module
    structure on $T_x\TM$ for which subspaces $T_x\FF$ and $T_x\FF^\perp$ are $\g$-submodules.
  \end{itemize}
 \end{proposition}

 We consider the $\g$-valued $1$-form $\theta$ on $\TM$ which is defined, at every $x\in\TM$, using the
 composition of the projection $T_x\TM\to T_x\FF$ and the isomorphism of $T_x\FF$ with $\g$. We also
 consider the $\g$-valued $2$-form given by $\Theta=d\theta|_{\wedge^2T\FF^\perp}$.

 The proof of the following result can be found in \cite[p. 239]{Q-TINC}.
 \begin{lemma}\label{lem: integrable}
  Let $G$, $M$ and $S$ be as in Proposition \ref{prop: g(x)}. If we assume that the $G$-orbits are non-degenerate,
  then:
  \begin{itemize}
   \item[$(1)$] For every $x\in S$, the maps $\theta_x:T_x\TM\to\g$ and $\Theta_x:\wedge^2T_x\FF^\perp\to\g$ are
    both homomorphism of $\g$-modules, for the $\g$-module structures from Proposition \ref{prop: g(x)}.
   \item[$(2)$] The normal bundle $T\FF^\perp$ is integrable if and only if $\Theta=0$.
  \end{itemize}
 \end{lemma}

 \begin{remark}
  Let $G$, $M$ and $S$ be as in Proposition \ref{prop: g(x)}, suppose that the $G$-orbits on $M$ are
  non-degenerate. From the precedent Lemma and the analyticity of the elements involved in this work we have
  two possible cases:
  \emph{(i) $\Theta\equiv0$ and then $T\FF^\perp$ is integrable}, or \emph{(ii) $\Theta\neq0$ for almost $x\in\TM$}.
 \end{remark}

 In what it follows we will assume that the $G$-orbits are non-degenerate with respect to the
 pseudo-Riemannian metric. Hence, the $\TG$-orbits on $\TM$ are non-degenerate as well and we have a direct sum
 decomposition $T\TM=T\FF\oplus T\FF^\perp$ . The non-degeneracy of the orbits is ensured for manifolds of low
 dimension with respect to the dimension of the Lie group acting on this, that is the result of the
 following Lemma that can be founded in \cite{Q-TINC}.

 \begin{lemma}[{\cite[Lemma 2.7]{Q-TINC}}]\label{lem: non-degen}
  Let $G$ be a connected non-compact simple Lie group acting isometrically and with a dense orbit on a connected
  finite volume pseudo-Riemannian manifold $M$. If $\dim(M)<2\dim(G)$, then the bundles $T\FF$ and $T\FF^\perp$
  have fibers that are non-degenerated with respect to the metric on $M$.
 \end{lemma}

 For the $G$-action as in Proposition \ref{prop: g(x)}, we consider $\TM$ endowed with the $\TG$-action obtained
 by lifting the $G$-action on $M$. Let us denote by $\HH$ the Lie subalgebra of $\Kill(\TM)$ consisting of the
 fields that centralize the $\TG$-action on $\TM$. Our first lemma involving $\HH$ is about an embedding of the
 Lie algebra $\g$ into $\HH$. Such result %of the following Lemma (proved in \cite{OQ}) which
 allows us to apply representation theory to the study of $\HH$.
 \begin{lemma}[{\cite[Lemma 1.7]{OQ}}]\label{lem: embed}
  Let $S$ as in Proposition \ref{prop: g(x)}. Then, for every $x\in S$ and for $\rho_x$ given as in Proposition
  \ref{prop: g(x)}, the map $\widehat{\rho}_x:\HH\to\Kill(\TM)$ defined as:
  \begin{equation*}
   \widehat{\rho}_x(X)=\rho_x(X)+X^*,
  \end{equation*}
  is an injective homomorphism of Lie algebras whose image $\GG(x)$ lies in $\HH$. In particular, $\widehat{\rho}_x$
  induces on $\HH$ a $\g$-module structure such that $\GG(x)$ is a submodule isomorphic to $\g$.
 \end{lemma}
 \begin{proof}
  First, observe that by the identity in Proposition \ref{prop: g(x)}(3), one can see that the image of
  $\widehat{\rho}_x$ lies in $\HH$.

  To prove that $\widehat{\rho}_x$ is a homomorphism of Lie algebras we apply Proposition \ref{prop: g(x)} as
  follows: for $X,Y\in\g$ we have
  \begin{eqnarray*}
   [\widehat{\rho}_x(X),\widehat{\rho}_x(Y)]&=&[\rho_x(X)+X^*,\rho_x(Y)+Y^*]\\
    &=&[\rho_x(X),\rho_x(Y)]+[X,Y]^*+[X,Y]^*+[X^*,Y^*]\\
    &=&\rho_x([X,Y])+[x,Y]^*\\
    &=&\widehat{\rho}_x([X,Y]).
  \end{eqnarray*}

  From the definition of $\widehat{\rho}_x$ we observe that $\widehat{\rho}_x(X)=0$ implies that
  $X^*=0$, which in turns yields $X=0$, because the $G$-action is locally free. Therefore, the last claim
  of our lemma is now clear.
 \end{proof}

 The following lemma relates the structure of $\g$-module of $\HH$ to that of $T_x\TM$.

 \begin{lemma}\label{lem: evaluation}
  Let $S$ as in Proposition \ref{prop: g(x)}. Consider $T_x\TM$ and $\HH$ endowed with the $\g$-module structure
  given by Proposition \ref{prop: g(x)}(4) and Lemma \ref{lem: embed}, respectively. Then, for almost every $x\in S$,
  the evaluation map:
  \begin{equation*}
   ev_x:\HH\to T_x\TM, \qquad Z\mapsto Z_x
  \end{equation*}
  is a homomorphism of $\g$-modules that satisfies $ev_x(\GG(x))=T_x\FF$. Furthermore, for almost every $x\in S$ we
  have $ev_x(\HH)=T_x\TM$.
 \end{lemma}
 \begin{proof}
  For every $x\in S$, if we let $Z\in\HH$ and $X\in\g$ be given, then:
  \begin{equation*}
   ev_x(X\cdot Z)=[\widehat{\rho}_x(X),Z]_x=[\rho_x(X)+X^*,Z]_x=[\rho_x(X),Z]_x=X\cdot Z_x=X\cdot ev_x(Z)
  \end{equation*}
  where we have used Lemma \ref{lem: homom to so} and the definition of the $\g$-module structures involved, thus
  proving the first part. The second part is proved by Lemma $4.1$ in \cite{Zimmer1}, using the transitivity of
  $\h$  on an open dense conull set in $M$.
 \end{proof}

 The next lemma shows the existence of a relation between isometries and Killing fields for complete
 manifolds. It also shows that every Lie algebra containing Killing vector fields can be obtained from an
 isometric right action, and whose proof can be found in \cite{OQ}.

 \begin{lemma}\label{lem: right act}
  Let $N$ be a complete pseudo-Riemannian manifold and $H$ a simply connected Lie group with Lie algebra
  $\h$. If $\psi:\h\to\Kill(N)$ is a homomorphism of Lie algebras, then there exists an isometric right
  $H$-action, $N\times H\to N$, such that $\psi(X)=X^*$, for every $X\in\h$. Furthermore, if $N$ is analytic,
  then the $H$-action is analytic as well.
 \end{lemma}

\section{Analysis and properties of isometric \texorpdfstring{$\SL(3,\R)$}-actions}\label{sect: Centralizer II}

 In this section we assume that $G=\SL(3,\R)$, which is a connected, non-compact simple Lie group, that the
 dimension of $M$ satisfies that $8<\dim(M)\leq14$ and we remain the
 other hypotheses of Proposition \ref{prop: g(x)}. Hence, $\g=\sli(3,\R)$ and, from Lemma \ref{lem: non-degen},
 we have that the orbits of the action on the manifold $M$ are non-degenerated.

 With our hypotheses in the previous paragraph we have the following result.

 \begin{theorem}\label{lem: integrability}
  With the hypotheses in Theorem A, if we assume that our normal bundle $T\FF^\perp$
  is integrable then $\TM\cong\TSL(3,\R)\times\TN$, where $\TN$ is a complete pseudo-Riemannian manifold.
 \end{theorem}
 \begin{proof}
  This is a direct consequence of the main results in \cite{Q-TINC}.
 \end{proof}

 The previous theorem induces to analyze the case when the normal bundle $T\FF^\perp$ is not integrable. Hence,
 from now on we assume that the normal bundle to the foliation is not integrable.

 Using the analyticity of our hypotheses we have the following result.

 \begin{lemma}\label{lem: dim14}
  Let $S$ as in Proposition \ref{prop: g(x)}. Let $x\in\sS$ be, consider $T_x\FF^\perp$ endowed with
  the $\sli(3,\R)$-module structure given by Proposition \ref{prop: g(x)}(4). Then, for almost every
  $x\in\sS$, $T_x\FF^\perp$ is isomorphic to $\R^3\oplus\R^{3*}$ and $\dim(M)=14$. In particular, the
  algebra $\so(T_x\FF^\perp)$ is isomorphic to $\so(3,3)$ as a Lie algebra and as a $\sli(3,\R)$-module.
 \end{lemma}
 \begin{proof}
  By the non-integrability of $T\FF^\perp$ then, from Lemma \ref{lem: integrable}(2), we have that the $2$-form
  $\Theta$ is not equal to zero. Since this $2$-form is analytic, then it vanishes on a proper
  analytic subset of null measure. Hence, $\Theta_x\neq0$ for almost every $x\in S$.

  Choose and fix an element $x\in S$ such that $\Theta_x\neq0$.
  By the definition of $\Theta_x$ we have that $T_x\FF^\perp$ is a non-trivial $\sli(3,\R)$-module
  preserving a non-degenerate symmetric bilinear form
  with $1<\dim(T_x\FF^\perp)\leq6$, then from Lemma \ref{lem: Bil Form} and Proposition \ref{prop: g(x)}(4),
  $T_x\FF^\perp\simeq\R^3\oplus\R^{3*}$ as a $\sli(3,\R)$-module. Hence,
  $\dim(T_x\FF^\perp)=6$ and, therefore, $\dim(M)=14$.

  On the other hand, by Lemma \ref{lem: Bil Form}, the representation of $\sli(3,\R)$ on $T_x\FF^\perp$
  defines a non-trivial homomorphism of Lie algebras $\so(T_x\FF^\perp)\to\so(3,3)$, which is also an
  isomorphism of $\sli(3,\R)$-modules. Since $\so(3,3)$ is
  a simple Lie algebra, this latter homomorphism is injective and so is an isomorphism.
 \end{proof}

 The results of the previous lemma allow us to obtain a decomposition of the centralizer $\HH$ of the
 $\TSL(3,\R)$-action into submodules related to the pseudo-Riemannian metric structure on $\TM$. First, recall
 that Lemma \ref{lem: embed} induces on $\HH$ a structure of $\sli(3,\R)$-module. % structure. % on $\HH$.

 \begin{lemma}\label{lem: decomposition}
  Let $S$ be as in Proposition \ref{prop: g(x)}. Then, for almost every $x\in S$ there is a decomposition of $\HH$
  into $\sli(3,\R)$-submodules, $\HH=\GG(x)\oplus\HH_0(x)\oplus\WW(x)$, satisfying:
  \begin{itemize}
   \item[1)] $\GG(x)=\widehat{\rho}_x(\sli(3,\R))$ is a Lie subalgebra of $\HH$ isomorphic to $\sli(3,\R)$ and
    $ev_x(\GG(x))=T_x\FF$.
   \item[2)] $\HH_0(x)=\ker(ev_x)$ is a Lie subalgebra and a $\sli(3,\R)$-module of $\HH$, isomorphic to a
    subset of $\so(T_x\FF^\perp)$.
   \item[3)] $\WW(x)$ is isomorphic to $\R^3\oplus\R^{3*}$ as $\sli(3,\R)$-module and $ev_x(\WW(x))=T_x\FF^\perp$.
  \end{itemize}
  In particular, the evaluation map (at $x$) defines an isomorphism of $\sli(3,\R)$-modules
  $\GG(x)\oplus\WW(x)\to T_x\FF\oplus T_x\FF^\perp$, preserving the summands in that order. In 2), we have
  that the considered isomorphism is as a Lie algebra and as an $\sli(3,\R)$-module.
 \end{lemma}
 \begin{proof}
  Let us choose and fix an element $x\in S$ that satisfies Lemma \ref{lem: evaluation} and Lemma \ref{lem: dim14}.
  By Lemma \ref{lem: embed} we conclude that $\GG(x)=\widehat{\rho}_x(\sli(3,\R))$ is a Lie algebra isomorphic
  to $\sli(3,\R)$.

  Define $\HH_0(x)=\ker(ev_x)$. By Lemma \ref{lem: evaluation}, it follows that $\HH_0(x)$ is an
  $\sli(3,\R)$-submodule of $\HH$. On the other hand, since $\HH_0(x)=\HH\cap\Kill_0(\TM,x)$, one can see that
  $\HH_0(x)$ is a Lie subalgebra of $\HH$.

  Since the elements of $\GG(x)$ are of the form $\rho_x(X)+X^*$, with $X\in\sli(3,\R)$.
  Hence, for any such element we have $ev_x(\rho_x(X)+X^*)=X_x^*$. That and the condition
  $ev_x(\rho_x(X)+X^*)=0$ imply that $X=0$.
  In other words, $\GG(x)\cap\HH_0(x)=0$; therefore there exists an $\sli(3,\R)$-submodule complementary
  $\WW'(x)$ to $\GG(x)\oplus\HH_0(x)$ in $\HH$. Note, we have an isomorphism from $\GG(x)\oplus\WW'(x)$
  onto $T_x\TM$ via the evaluation map. We choose $\WW(x)$ as the inverse image of $T_x\FF^\perp$
  under our previous isomorphism. We have our desired decomposition of $\HH$ into $\sli(3,\R)$-submodules.

  Let $\Kill_0(\TM,x,\FF)$ be the Lie algebra of Killing vector fields on $\TM$ which preserves the foliation $\FF$
  and vanish at $x$. Note that every vector field in $\Kill_0(\TM,x,\FF)$ leaves invariant the normal bundle.
  On the other hand, the map $\lambda_x$ restricted to $\Kill_0(\TM,x,\FF)$ induces the following
  homomorphism of Lie algebras:
  \begin{equation*}
   \lambda_x^\perp:\Kill_0(\TM,x,\FF)\to\so(T_x\FF^\perp), \qquad X\mapsto\lambda_x(X)|_{T_x\FF^\perp}.
  \end{equation*}
  Observe that both $\rho_x(\sli(3,\R))$ and $\HH_0(x)$ lie inside of $\Kill_0(\TM,x,\FF)$.

  \emph{Claim $1$: $\lambda_x^\perp$ is injective when is restricted to $\sli(3,\R)(x)$}. By our choice of the
  element $x\in S$ and the results in Lemma \ref{lem: dim14}, the map
  $\lambda_x^\perp\circ\rho_x:\sli(3,\R)\to\so(T_x\FF^\perp)$ is a non-trivial homomorphism of Lie algebras.
  Since $\sli(3,\R)$ is a simple Lie algebra then the map $\lambda_x^\perp$ restricted to $\sli(3,\R)$ is
  injective.

  \emph{Claim $2$: $\lambda_x^\perp$ restricted to $\HH_0(x)$ is injective}. Recall that pseudo-Riemannian metric
  structures are $1$-rigid structures (see \cite{GCT-CQ}), therefore a Killing vector space is completely
  determined by its $1$-jet at $x$. If $Z\in\HH_0(x)$ is given, then $Z_x=ev_x(Z)=0$; so it is determined
  by the values $[Z,V]_x$, for $V$ vector field on a neighborhood of $x$. Since $Z$ is in the centralizer of the
  $\TSL(3,\R)$-action, then $[Z,X^*]_x=0$ for all $X\in\sli(3,\R)$, so $[Z,V]_x=0$ when $V_x\in T_x\FF$. Hence,
  if $[Z,V]_x=0$ when $V_x\in T_x\FF^\perp$, this implies that $Z=0$. Therefore, we have that $\lambda_x^\perp$
  is injective when it is restricted to $\HH_0(x)$.

  On the other hand, if $X\in\sli(3,\R)$ and $Y\in\HH_0(x)$ then
  \begin{eqnarray*}
   \lambda_x^\perp(X\cdot Y)&=&\lambda_x^\perp([\widehat{\rho}_x(X),Y])=\lambda_x^\perp([\rho_x(X)+X^*,Y])\\
    &=&\lambda_x^\perp([\rho_x(X),Y])=[\lambda_x^\perp(\rho_x(X)),\lambda_x^\perp(Y)]\\
    &=&X\cdot \lambda_x^\perp(Y).
  \end{eqnarray*}
  That shows that the map $\lambda_x^\perp$ restricted to $\HH_0(x)$ is a homomorphism of $\sli(3,\R)$-modules.
 \end{proof}

 \begin{remark}
  By Lemma \ref{lem: decomposition} we have that $\HH_0(x)$ is a subalgebra, and an $\sli(3,\R)$-submodule,
  isomorphic to $\lambda_x^\perp\big(\HH_0(x)\big)\subset\so(T_x\FF^\perp)$. On the other hand, since
  $\so(T_x^\perp)$ is isomorphic to $\so(3,3)$ (and hence to $\sli(4,\R)$), then $\HH_0(x)$ is isomorphic
  to one of the Lie subalgebras in Corollary \ref{cor: subs}.
 \end{remark}

 By Lemma A.5 in \cite{OQ-Upq}, we have that $\wedge^2T_x\FF^\perp$ is isomorphic to
 $\so(T_x\FF^\perp)$ as $\so(T_x\FF^\perp)$-module. Thus, from the definition of the map $\Theta_x$ in
 Lemma \ref{lem: integrable} and Lemma A.5 in \cite{OQ-Upq}, we can consider $\Theta_x$ as a map from
 $\so(T_x\FF^\perp)$ to $\sli(3,\R)$.

 More properties about the subalgebra $\HH_0(x)$ can be obtained considering the map $\Theta_x$ as in the
 previous paragraph, one of these is contained in the following result which appears in Proposition $3.10$
 and Proposition $3.11$ in \cite{OQ-Upq}.

 \begin{theorem}\label{th: H_0}
   For almost every $x\in\sS$ as in Lemma \ref{lem: decomposition}. $\lambda^\perp_x(\HH_0(x))$ is a
   $\lambda^\perp_x(\GG(x))$-submodule and a Lie algebra of $\so(T_x\FF^\perp)$ that satisfies
   \begin{equation*}
     [\lambda^\perp_x(\HH_0(x)),\so(T_x\FF^\perp)]\subset\ker(\Theta_x).
   \end{equation*}
 \end{theorem}
 \begin{proof}
   Let $x\in\sS$ be an element satisfying Lemma \ref{lem: decomposition}, so it does
   Lemma \ref{lem: evaluation} and Lemma \ref{lem: dim14}.

   The element $x\in\sS$ has the required properties of the hypotheses of Proposition $3.10$ in
   \cite{OQ-Upq}. Therefore, the proof of this theorem is similar to that of the above mentioned
   proposition.
 \end{proof}

 \begin{remark}\label{rem: 0 or R}
   Decomposing $\so(T_x\FF^\perp)$ (isomorphic to $\so(3,3)$) as a direct sum of irreducible
   $\sli(3,\R)$-submodules, and its corresponding bracket product, Lemma \ref{lem: D irred} yields
   that the only possibilities for $\HH_0(x)$ satisfying Theorem \ref{th: H_0} are $0$, or isomorphic to $\R$.
 \end{remark}

\section{Structure of the Centralizer and its consequences}\label{sect: Structure}

 The previous section was devoted to show all the possible values that $\HH_0(x)$ can take. This section
 we analyze the implications of all these possible cases. Here, we assume the same hypotheses and notation
 of Lemma \ref{lem: decomposition}.
 
 First, fixed an arbitrarily element $x\in\sS$ as in Lemma \ref{lem: decomposition} and Theorem \ref{th: H_0},
 guaranty
 (by Lemma \ref{lem: dim14} and Lemma \ref{lem: decomposition}) we can choose subspaces $\VV(x)$ and $\VV^*(x)$
 of $\WW(x)$ such that $\WW(x)=\VV(x)\oplus\VV^*(x)$ and $\GG(x)$ acts on $\VV(x)$ (resp. $\VV^*(x)$) as
 $\sli(3,\R)$ acts on $\R^3$ (resp. $\R^{3*}$).

 Thus, $\HH$ is an $\sli(3,\R)$-module and on properties of the evaluation map in Lemma
 \ref{lem: decomposition}, we have the following properties:
 \begin{align}
  [\GG(x),\HH_0(x)] &\subseteq \HH_0(x) \label{eq: [G,H0]cH0}\\
  [\GG(x),\WW(x)] &= \WW(x) \label{eq: [G,W]=W}\\
  [\HH_0(x),\WW(x)] &\subseteq \HH_0(x)\oplus\WW(x) \label{eq: [H0,W]cH0+W}\\
  [\HH_0(x),\HH_0(x)] &\subseteq \HH_0(x) \label{eq: [H0,H0]cH0}
 \end{align}

 In particular, when $\HH_0(x)$ is isomorphic to $\R$ or to $0$, the contention \eqref{eq: [H0,W]cH0+W}
 is strict. That is a consequence of the following lemma.
 \begin{lemma}\label{lem: equalities}
  As in Lemma \ref{lem: decomposition}, if $\HH_0(x)$ is isomorphic either to $\R$, $\sli(3,\R)$,
  $\sli(3,\R)\oplus\R$ or $\sli(4,\R)$ then $[\HH_0(x),\WW(x)]=\WW(x)$.
 \end{lemma}

 Next, we analyze the different possibilities of $\HH_0(x)$ satisfying Theorem \ref{th: H_0}. Recall
 that $\HH_0(x)$ is isomorphic to one of the subalgebras of Remark \ref{rem: 0 or R}.

\subsubsection{$\HH_0(x)=0$}\label{ss: subsub H0=0}
 \begin{lemma}\label{lem: H0=0}
  Let $\sS$ be as in Lemma \ref{lem: decomposition}. If $x\in\sS$ and $\HH_0(x)=0$ then, one of the following
  occurs:
  \begin{enumerate}
   \item The Radical of $\HH$ is $\WW(x)=\VV(x)\oplus\VV^*(x)$.
   \item $\HH=\GG(x)\oplus\VV(x)\oplus\VV^*(x)$ is isomorphic to $\g_{2(2)}$.
  \end{enumerate}
 \end{lemma}
 \begin{proof}
  Let us choose an arbitrary but fixed element $x\in\sS$ such that $\HH_0(x)=0$, as in Lemma
  \ref{lem: decomposition}. Hence, $\HH=\GG(x)\oplus\WW(x)=\GG(x)\oplus\VV(x)\oplus\VV^*(x)$.

  Since $\GG(x)$ is a simple Lie algebra, we can choose $\s$ a Levi factor of $\HH$ which contains $\GG(x)$.
  Since $\GG(x)\subseteq\s$ and considering the structure of $\HH$ as an $\sli(3,\R)$-module, then $\s$ is
  also an $\sli(3,\R)$-module.

  Let $W$ be an $\sli(3,\R)$-module of $\HH$ such that $\s=\GG(x)\oplus W$. Moreover,
  since $\rad(\HH)$ is an ideal of $\HH$ this induces the decomposition of $\HH$ as a direct sum of
  $\sli(3,\R)$-modules:
  \begin{equation*}
    \HH=\s\oplus\rad(\HH)=\GG(x)\oplus W\oplus\rad(\HH)
  \end{equation*}
  that we compare with its decomposition into irreducible $\sli(3,\R)$-submodules from Lemma
  \ref{lem: decomposition}
  \begin{equation*}
   \HH=\GG(x)\oplus\VV(x)\oplus\VV^*(x).
  \end{equation*}

  The properties of representations of Lie algebras and the decompositions of $\HH$ imply that one of
  the following must occur:
  \begin{itemize}
   \item[(a)] $\s=\GG(x)\oplus\VV(x)$ and $\rad(\HH)=\VV^*(x)$.
   \item[(b)] $\s=\GG(x)\oplus\VV^*(x)$ and $\rad(\HH)=\VV(x)$.
   \item[(c)] $\s=\GG(x)$ and $\rad(\HH)=\VV(x)\oplus\VV^*(x)$.
   \item[(d)] $\HH=\GG(x)\oplus\VV(x)\oplus\VV^*(x)$ is semisimple.
  \end{itemize}
 
  \textbf{Suppose the case (a) is satisfied:}
  \begin{equation*}
   \s=\GG(x)\oplus\VV(x) \quad \text{and} \quad\rad(\HH)=\VV^*(x).
  \end{equation*}

  Because $\s=\GG(x)\oplus\VV(x)$ is a semisimple Lie algebra, then $\s$ is a finite direct pro\-duct,
  of simple ideals $\h_1\times\h_2\times\dots\times\h_k$. Since, every ideal is invariant
  by $\GG(x)$ then these ideals are $\sli(3,\R)$-modules. Representation properties of $\sli(3,\R)$ and
  the decomposition of $\s$ in a direct sum of irreducible $\sli(3,\R)$-modules yield that $k\leq2$.

  If $k=2$, $\s=\h_1\times\h_2$. By \eqref{eq: [G,W]=W} we have that $[\GG(x),\VV(x)]\subseteq\VV(x)$.
  Additionally, since $\VV(x)$ is isomorphic to $\R^3$ as $\sli(3,\R)$-module then, $[\VV(x),\VV(x)]\neq0$
  implies that $[\VV(x),\VV(x)]\simeq{\R^3}^*$. Which shows that $[\VV(x),\VV(x)]=0$.
  Therefore $\VV(x)$ is an ideal of
  $\s=\GG(x)\oplus\VV(x)$. Without loss of generality, we assume $\h_2=\VV(x)$, then the simple ideal $\h_2$
  must be an abelian ideal, thus $k=1$.

  However, if $k=1$, $\s$ is a simple Lie algebra. Thence, $\s=\GG(x)\oplus\VV(x)$ is a $11$-dimensional real simple
  Lie algebra. Therefore, $\s^\C$ is a complex simple Lie algebra with complex dimension $11$ and
  \cite[p. 516]{Helgason}, showed it cannot happen. Then, we have proved that case (a) cannot happen.

  \textbf{Assume case (b) is satisfied:}
  \begin{equation*}
   \s=\GG(x)\oplus\VV^*(x) \quad \text{and} \quad\rad(\HH)=\VV(x).
  \end{equation*}

  This case is ruled out by a similar argument of case (a).

  \textbf{Suppose case (c) is satisfied:}
  \begin{equation*}
   \s=\GG(x) \quad \text{and} \quad\rad(\HH)=\VV(x)\oplus\VV^*(x).
  \end{equation*}

  Since $\VV(x)\simeq\R^3$ and $\VV^*(x)\simeq\R^{3*}$ as $\sli(3,\R)$-modules then by properties of representation
  of $\sli(3,\R)$ we have that $[\VV(x),\VV^*(x)]$ is isomorphic to a subspace of $\sli(3,\R)\oplus\R$ as
  a $\sli(3,\R)$-module.
  Therefore, we have that $[\VV(x),\VV(x)]\subseteq\VV^*(x)$, $[\VV(x),\VV^*(x)]=0$ and
  $[\VV^*(x),\VV^*(x)]\subseteq\VV(x)$.
  And by the solvability of $\rad(\HH)$ we have that $[\VV(x),\VV(x)]=0$ or $[\VV^*(x),\VV^*(x)]=0$.

  \textbf{If case (d) is satisfied:}
  \begin{equation*}
   \HH=\GG(x)\oplus\VV(x)\oplus\VV^*(x) \quad \text{is a simple Lie algebra}.
  \end{equation*}

  Using the same argument as in case (a), $\HH$ is a direct product of a finite number
  of simple ideals $\h_1\times\h_2\times\dots\times\h_k$ where every ideal is an $\sli(3,\R)$-module and $k\leq3$.

  If $k=3$, $\HH=\h_1\times\h_2\times\h_3$, we can assume, reindexing if necessary, that $\h_3=\VV^*(x)$ and
  $\h_1\times\h_2=\GG(x)\oplus\VV(x)$. Then $[\VV^*(x),\VV^*(x)]\subseteq\h_3$, but similar to case (a) we have that
  $[\VV^*(x),\VV^*(x)]=0$, so $\h_3$ is an abelian Lie algebra and this is not possible. Therefore $k\leq2$.

  If $k=2$, $\HH=\h_1\times\h_2$, after decomposing $\h_1$ and $\h_2$ as the direct sum of $\sli(3,\R)$-modules,
  and reindexing if necessary, we can assume that $\h_1$ is an irreducible $\sli(3,\R)$-module and
  $\h_2=V_1\oplus V_2$, where $V_1$ and $V_2$ are irreducible $\sli(3,\R)$-modules.
  We can also assume that $\VV^*(x)\subset\h_2$ and $\VV(x)\oplus\VV^*(x)=\h_2$, because as in the previous
  case, it cannot happen that $\h_1=\VV(x)$. Then $\GG(x)\subseteq\h_1$ and $[\GG(x),\WW(x)]=[\h_1,\h_2]=0$ that
  contradicts the equation \eqref{eq: [G,W]=W}.

  If $k=1$, $\HH=\GG(x)\oplus\VV(x)\oplus\VV^*(x)$ is a real simple Lie algebra of dimension $14$. Then, $\HH$ is
  the realification of a complex simple Lie algebra of dimension $7$ or its complexification, $\HH^\C$, is a
  complex simple Lie algebra. Since, by \cite[p. 516]{Helgason}, there is not a complex simple Lie algebra
  of dimension $7$. So, $\HH^\C$ is a complex simple Lie algebra with $\dim_\C(\HH^\C)=14$. Then,
  $\HH^\C\cong\g_2$. From here, $\HH$ is isomorphic to a real form of $\g_2$.

  On the other hand, we recall that $\HH$ contains a Lie subalgebra isomorphic to $\sli(3,\R)$, that is simple
  and non-compact. Then, $\HH$ is non-compact. Otherwise, exercise $4(ii)$ in the page $152$ of \cite{Helgason},
  would imply that $\sli(3,\R)$ is compact, which is a contradiction.

  Since, by \cite[p. 518]{Helgason}, there is only a non-compact real form of $\g_2$, namely $\g_{2(2)}$. Then
  \begin{equation*}
   \HH=\GG(x)\oplus\VV(x)\oplus\VV^*(x)\simeq\g_{2(2)}.
  \end{equation*}
 \end{proof}

\subsubsection{$\HH_0(x)\simeq\R$}
 \begin{lemma}\label{lem: H0=R}
  Let $\sS$ be as in Lemma \ref{lem: decomposition}. Let $x\in\sS$ be such that $\HH_0(x)$ is isomorphic to
  $\R$ as an $\sli(3,\R)$-module then one of the following occurs
  %If $\HH_0(x)$ is isomorphic to $\R$ as $\sli(3,\R)$-module for some $x\in\sS$, one of the following occurs:
 \begin{enumerate}
  \item[(1)] The radical of $\HH$ is $\HH_0(x)\oplus\VV(x)\oplus\VV^*(x)$ where $\VV(x)\oplus\VV^*(x)$ is a Lie
   subalgebra.
  \item[(2)] $\HH=\GG(x)\oplus\HH_0(x)\oplus\VV(x)\oplus\VV^*(x)$ is a simple Lie algebra isomorphic to $\sli(4,\R)$.
 \end{enumerate}
 \end{lemma}
 \begin{proof}
  Let us choose an arbitrary but fixed element $x\in\sS$, as in Lemma \ref{lem: decomposition}, such that
  $\HH_0(x)\simeq\R$.

  Choose $\s$ a Levi factor of $\HH$ that contains $\GG(x)$. Similar to the case $\HH_0(x)=0$, $\s$ is a
  $\sli(3,\R)$-submodule of $\HH$.

  Let $W$ be a $\sli(3,\R)$-submodule of $\HH$ such that $\s=\GG(x)\oplus W$.
  Since $\rad(\HH)$ is an ideal, this induces the next decomposition of $\HH$ as a direct sum of
  $\sli(3,\R)$-modules:
  \begin{equation*}
   \HH=\GG(x)\oplus W\oplus\rad(\HH)
  \end{equation*}
  that we compare with the decomposition of irreducible $\sli(3,\R)$-modules in Lemma \ref{lem: decomposition}
  \begin{equation*}
   \HH=\GG(x)\oplus\HH_0(x)\oplus\VV(x)\oplus\VV^*(x).
  \end{equation*} %from Lemma \ref{lem: decomposition}

  By the properties of representations of Lie algebras and the decomposition of
  $\HH=\GG(x)\oplus\HH_0(x)\oplus\VV(x)\oplus\VV^*(x)$, one of the following must occur:
  \begin{itemize}
   \item[(a)] $\s=\GG(x)\oplus\HH_0(x)\oplus\VV(x)$ and $\rad(\HH)=\VV^*(x)$.
   \item[(b)] $\s=\GG(x)\oplus\HH_0(x)\oplus\VV^*(x)$ and $\rad(\HH)=\VV(x)$.
   \item[(c)] $\s=\GG(x)\oplus\VV(x)\oplus\VV^*(x)$ and $\rad(\HH)=\HH_0(x)$.
   \item[(d)] $\s=\GG(x)\oplus\HH_0(x)$ and $\rad(\HH)=\VV(x)\oplus\VV^*(x)$.
   \item[(e)] $\s=\GG(x)\oplus\VV(x)$ and $\rad(\HH)=\HH_0(x)\oplus\VV^*(x)$.
   \item[(f)] $\s=\GG(x)\oplus\VV^*(x)$ and $\rad(\HH)=\HH_0(x)\oplus\VV(x)$.
   \item[(g)] $\s=\GG(x)$ and $\rad(\HH)=\HH_0(x)\oplus\VV(x)\oplus\VV^*(x)$.
   \item[(h)] $\HH=\GG(x)\oplus\HH_0(x)\oplus\VV(x)\oplus\VV^*(x)$ is semisimple.
  \end{itemize}

  \textbf{Suppose that case (a) is satisfied:}
  \begin{equation*}
   \s=\GG(x)\oplus\HH_0(x)\oplus\VV(x) \quad \text{and} \quad\rad(\HH)=\VV^*(x).
  \end{equation*}

  Recall that $\s=\GG(x)\oplus\HH_0(x)\oplus\VV(x)$ is a semisimple Lie algebra.

  Since $\VV(x)\simeq\R^3$ and $[\VV(x),\VV(x)]\subseteq\s$. If $[\VV(x),\VV(x)]\neq0$ then
  $[\VV(x),\VV(x)]$ is isomorphic to $\R^{3*}$. Therefore, the projection of $[\VV(x),\VV(x)]$ in
  $\GG(x)$, $\HH_0(x)$ and $\VV(x)$ is $0$. That implies that $[\VV(x),\VV(x) ]=0$. Hence, $\VV(x)$ is an abelian
  ideal of $\s$, which is a contradiction. So, case (a) cannot be possible.

  \textbf{Case (b) is not possible} and the proof is similar to (a).

  \textbf{Now suppose case (c) is satisfied:}
  \begin{equation*}
   \s=\GG(x)\oplus\VV(x)\oplus\VV^*(x) \quad \text{and} \quad\rad(\HH)=\HH_0(x).
  \end{equation*}

  Here, by Lemma \ref{lem: equalities}, we have that $[\HH_0(x),\VV(x)\oplus\VV^*(x)]=\VV(x)\oplus\VV^*(x)$.
  On the other hand, since $\rad(\HH)$ is an ideal of $\HH$, then $\VV(x)\oplus\VV^*(x)\subseteq\rad(\HH)$. So,
  this case cannot occur.

  \textbf{Suppose case (d) is satisfied:}
  \begin{equation*}
   \s=\GG(x)\oplus\HH_0(x) \quad \text{and} \quad\rad(\HH)=\VV(x)\oplus\VV^*(x).
  \end{equation*}

  From \eqref{eq: [G,H0]cH0} and \eqref{eq: [H0,H0]cH0}, we have that $\HH_0(x)$ is an abelian ideal of $\s$.
  Which is a contradiction. Then case (d) is not possible.

  \textbf{Suppose, the case (e) is satisfied:}
  \begin{equation*}
   \s=\GG(x)\oplus\VV(x) \quad \text{and} \quad\rad(\HH)=\HH_0(x)\oplus\VV^*(x).
  \end{equation*}

  Similar to the argument in case (a), $\VV(x)$ is an abelian ideal of $\s$. Therefore, this case is not possible.

  \textbf{Case (f)},
  \begin{equation*}
   \s=\GG(x)\oplus\VV^*(x) \quad \text{and} \quad\rad(\HH)=\HH_0(x)\oplus\VV(x),
  \end{equation*}
  cannot happen and the proof is similar to that of (e).

  \textbf{Now, we suppose (g) is satisfied:}
  \begin{equation*}
   \s=\GG(x) \quad \text{and} \quad\rad(\HH)=\HH_0(x)\oplus\VV(x)\oplus\VV^*(x).
  \end{equation*}

  By Lemma \ref{lem: equalities} we have $[\HH_0(x),\WW(x)]=\WW(x)$. Let $\pi_0:rad(\HH)\to\HH_0(x)$ be the
  projection map on the first component of $\rad(\HH)$, since
  $[\WW(x),\WW(x)]\subseteq\rad(\HH)$ then $\pi_0([\WW(x),\WW(x)])=0$.
  Otherwise $\rad(\HH)$ will not be solvable. In conclusion, $\VV(x)\oplus\VV^*(x)$ is a Lie subalgebra of $\rad(\HH)$.

  \textbf{If case (h) is satisfied:}
  \begin{equation*}
   \HH=\GG(x)\oplus\HH_0(x)\oplus\VV(x)\oplus\VV^*(x) \quad \text{is semisimple}.
  \end{equation*}

  Here, $\HH$ is a finite direct product $\h_1\times\h_2\times\dots\times\h_k$ of simple ideals which also are
  $\sli(3,\R)$-modules with $k\leq4$.

  If $k=4$, $\HH=\h_1\times\h_2\times\h_3\times\h_4$. Without loss of generality we can assume that
  $\h_4=\HH_0(x)$ is a simple ideal of $\HH$. Since
  $\HH_0(x)\simeq\R$ is abelian, this is a contradiction. Therefore $k=4$ cannot be possible and $k\leq3$.

  If $k=3$, $\HH=\h_1\times\h_2\times\h_3$. Suppose, reindexing if necessary, that $\h_1$ and $\h_2$ are
  irreducible $\sli(3,\R)$-modules and $\h_3=V_1\oplus V_2$, where $V_1$ and $V_2$ are irreducibles
  $\sli(3,\R)$-modules. We can also assume that $\HH_0(x)\subset\h_3$ then, by Lemma \ref{lem: equalities},
  $\VV(x)\oplus\VV^*(x)\subset\h_3$. That implies that $\h_1$ or $\h_2$ is equal to $0$. Therefore, this
  case is not possible and $k\leq2$.

  If $k=2$, $\HH=\h_1\times\h_2$. As in case $k=3$, assume $\HH_0(x)\oplus\WW(x)\subset\h_2$.
  On the other hand, by \eqref{eq: [G,W]=W}, $[\GG(x),\WW(x)]=\WW(x)$ then
  $\GG(x)\subset\h_2$ and $\h_1=0$, which is a contradiction. Hence, this case cannot be possible and $k=1$.

  If $k=1$, $\HH$ is a real simple Lie algebra of dimension
  $15$. Therefore, $\HH^\C$ is a complex simple Lie algebra with $\dim_\C(\HH^\C)=15$.
  Then, by \cite[p. 516]{Helgason}, $\HH^\C\simeq\sli(4,\C)$. So, $\HH$ is
  isomorphic to a non-compact real form of $\sli(4,\C)$.

  From \cite[Table V]{Helgason}, the only non-compact real forms of $\sli(4,\C)$ are $\su(1,3)$,
  $\su(2,2)$, $\su^*(4)$ and $\sli(4,\R)$.

  Then, $\HH$ is isomorphic to one of the previous Lie algebras. Recall that $\HH$ contains
  a simple Lie subalgebra isomorphic to $\sli(3,\R)$. From here $2=\Rank_\R(\sli(3,\R))\leq\Rank_\R(\HH)$.
  By \cite[Table V]{Helgason}, we have $\Rank_\R(\su(1,3))=\Rank_\R(\su^*(4))=1$.
  Then, $\HH$ cannot be isomorphic to either $\su(1,3)$ or $\su^*(4)$. On the other hand, page $519$ of
  \cite{Helgason} shows that $\su(2,2)\simeq\so(4,2)$. So, if $\HH\simeq\su(2,2)$ then $\sli(3,\R)$ is isomorphic
  to a Lie subalgebra of $\so(4,2)$. In this case $\sli(3,\R)$ would have a non-trivial representation on a
  $6$-dimensional vector space that preserves a non-degenerate symmetric bilinear form of signature $(4,2)$.
  By Lemma \ref{lem: Bil Form} this cannot be possible. Thus, $\HH\simeq\su(2,2)$ is not possible.
  Then $\HH\simeq\sli(4,\R)$.
 \end{proof}

\section{Proof of the Main Theorem}\label{sect: Proof}
 In this section we assume $M$ is a connected analytic pseudo-Riemannian manifold with $\dim(M)=14$ and finite
 volume. We also assume that $\SL(3,\R)$ acts isometric and analytically on $M$ with a dense orbit, therefore
 the action is locally free, such that the normal bundle to the foliation (obtained by the action) $T\FF^\perp$
 is not integrable. We study the
 structure of the manifold $M$ through the analysis of the different possibilities of $\HH$ obtained %in the results
 in Section \ref{sect: Structure}. As in the previous section, we use the notation of Lemma \ref{lem: decomposition}.

 Let $x\in\TM$ be an element that satisfies Lemma \ref{lem: decomposition} and Theorem \ref{th: H_0}.
 By Lemma \ref{lem: H0=0} and Lemma \ref{lem: H0=R} we can reduce the structure of $\HH$ to
 only $2$ options: (1) $\WW(x)\subseteq\rad(\HH)$ is a subalgebra or (2) $\HH$ is a simple Lie algebra.
 In this section we analyze these cases and, assuming that $T\FF^\perp$ is not integrable then, we will see that
 we can eliminate the first possibility.

\subsection{\texorpdfstring{$\WW(x)\subseteq\rad(\HH)$}  \textbf{ is a subalgebra.}}\label{subs: w nilp}

 Here, we assume that $\WW(x)$ is a subalgebra of $\rad(\HH)$ for some $x\in\TM$.

 Since $\rad(\HH)$ is a solvable Lie algebra then, as case (c) in Lemma \ref{lem: H0=0} and case (g) in
 \ref{lem: H0=R}, $\WW(x)=\VV(x)\oplus\VV^*(x)$ is a $2$-step nilpotent or an abelian Lie subalgebra.

 By Lemma \ref{lem: decomposition} and since $\WW(x)$ is a subalgebra,
 we have that $\GG(x)\oplus\WW(x)$ is a Lie subalgebra of $\HH$. Thus, $\GG(x)\oplus\WW(x)$ is isomorphic, as Lie
 algebra, to the semidirect product $\sli(3,\R)\ltimes\mathfrak{w}$, where $\mathfrak{w}$ is an $\sli(3,\R)$-module,
 $2$-step nilpotent or abelian Lie algebra isomorphic to $\WW(x)$. Choose an isomorphism of Lie algebras
 $\psi:\sli(3,\R)\ltimes\mathfrak{w}\to\HH(x)$ that maps $\sli(3,\R)$ onto $\GG(x)$ and $\mathfrak{w}$
 onto $\WW(x)$.

 Let $\TSL(3,\R)\ltimes\W$ be a simply connected Lie group such that
 Lie$(\TSL(3,\R)\ltimes\W)=\sli(3,\R)\ltimes\mathfrak{w}$, where the group structure on
 $\W$ is induced considering the action of $\sli(3,\R)$ on $\mathfrak{w}$. By Lemma \ref{lem: right act},
 there exists an analytic isometric right action of
 $\TSL(3,\R)\ltimes\W$ on $\TM$ such that $\psi(X)=X^*$ for every $X\in\sli(3,\R)\ltimes\mathfrak{w}$.

 Since $\HH$ centralizes the left $\TSL(3,\R)$-action, then the right
 $\big(\TSL(3,\R)\ltimes\W\big)$-action centralizes the left $\TSL(3,\R)$-action as well and preserves
 both $T\FF$ and $T\FF^\perp$.

 Using the right $\big(\TSL(3,\R)\ltimes\W\big)$-action on $\TM$, we consider the following map:
 \begin{equation}\label{eq: map from ra}
  p:\TSL(3,\R)\ltimes\W\to\TM, \qquad h\mapsto x\cdot h,
 \end{equation}
 for $h\in\TSL(3,\R)\ltimes\W$. This action is $\big(\TSL(3,\R)\ltimes\W\big)$-equivariant by the right
 action on its domain. If $e$ and $0$ denote the identity element in the subgroups $\TSL(3,\R)$ and $\W$,
 respectively, then
 \begin{eqnarray}\label{eq: diff of p}
  dp_{(e,0)}:&\sli(3,\R)\ltimes\mathfrak{w}\to\GG(x)\oplus\WW(x)\to T_x\TM&\\
  &X\mapsto X^* \mapsto X^*_x.& \nonumber
 \end{eqnarray}
 Since $\psi(X)=X^*$ for all $X\in\sli(3,\R)\ltimes\mathfrak{w}$, by Lemma \ref{lem: decomposition},
 $dp_{(e,0)}$ maps $\sli(3,\R)$ onto $T_x\FF$ and $\mathfrak{w}$ onto $T_x\FF^\perp$. Therefore, $p$ is a
 local diffeomorphism at $(e,0)$.

 For every $w\in\W$, let $R_w$ denote the map on $\TSL(3,\R)\ltimes\W$ and on $\TM$ given by the correspondence
 $y\mapsto y\cdot(e,w)$. Since $\W$ is a subgroup of $\TSL(3,\R)\ltimes\W$ we have that $R_w(\W)=\W$.

 Let $P=p(e\times\W)$, which defines a submanifold of $\TM$ in a neighborhood of $x=p(e,0)$. Here, by the
 previous remarks, we have that
 \begin{equation*}
   T_{p(e,0)}P=dp_{(e,0)}(T_{(e,0)}(e\times\W))=T_{p(e,0)}\FF^\perp,
 \end{equation*}
 which with the equivariance of $p$ implies that
 \begin{eqnarray*}
  T_{p(e,w)}P&=&dp_{(e,w)}(T_{(e,w)}(e\times\W))=dp_{(e,w)}(d(R_w)_{(e,0)}(T_{(e,0)}(e\times\W)))\\
   &=&d(R_w\circ p)_{(e,0)}(T_{(e,0)}(e\times\W))=d(R_w)_{p(e,0)}(T_{p(e,0)}P)\\
   &=&d(R_w)_{p(e,0)}T_{p(e,0)}\FF^\perp=T_{R_w(p(e,0))}\FF^\perp=T_{p(e,w)}\FF^\perp,
 \end{eqnarray*}
 we have used in the previous identities that $R_w$ preserves the bundle $T\FF^\perp$. This proves that
 $P$ is an integral submanifold of $T\FF^\perp$ passing through the element $x=p(e,0)$.

 By the left $\TSL(3,\R)$-action on $\TM$ we obtain by restriction to $P$ the following map:
 \begin{equation*}
  \phi:\TSL(3,\R)\times P\to\TM, \qquad (g,y)\mapsto g\cdot y,
 \end{equation*}
 whose differential at $(e,x)$ is given by: $X+v\mapsto X^*_x+v$, with $X\in\sli(3,\R)$ and $v\in T_xP$.
 This shows that the differential at $(e,x)$ is an isomorphism and therefore the map $\phi$ is a
 diffeomorphism from a neighborhood of $(e,x)$ onto
 a neighborhood of $x$. Since the left $\TSL(3,\R)$-action preserves both $T\FF$ and $T\FF^\perp$, there is an
 integral submanifold of $T\FF^\perp$ passing through every point in a neighborhood of $x$ in $\TM$. Thus, the tensor
 $\Theta$ considered in Lemma \ref{lem: integrable} vanishes in a neighborhood of $x$. Since all of our objects are
 analytic, this implies that $\Theta$ vanishes everywhere therefore Lemma \ref{lem: integrable} implies the
 integrability of $T\FF^\perp$ everywhere in $\TM$.

 This last conclusion contradicts the assumption about the integrability of $T\FF^\perp$.

\subsection{\texorpdfstring{$\HH$}  \textbf{ is a simple Lie algebra.}}\label{subs: H simple}

 Here, we assume $\HH(x)$ is a simple Lie algebra. By Lemma \ref{lem: H0=0} and Lemma \ref{lem: H0=R} we have
 two possibilities: (a) $\HH_0(x)=0$ and $\HH(x)\simeq\g_{2(2)}$ or (b) $\HH_0(x)\simeq\R$ and
 $\HH(x)\simeq\sli(4,\R)$

\subsubsection{$\HH\simeq\g_{2(2)}$}
 If $\HH_0(x)=0$ and $\HH$ is a simple Lie algebra then we have proved that $\HH\simeq\g_{2(2)}$, therefore
 our the following result:
 \begin{lemma}\label{lem: iso to g22}
  There is an isomorphism
  $$\psi:\g_{2(2)}=\sli(3,\R)\oplus\R^3\oplus\R^{3*}\to\GG(x)\oplus\VV(x)\oplus\VV^*=\HH$$
  of Lie algebras. In particular, we have that $\psi$ is an isomorphism
  of $\sli(3,\R)$-modules which preserves the summands in that order.
 \end{lemma}
 \begin{proof}
  Let $\psi:\g_{2(2)}\to\HH$ be an isomorphism of simple Lie algebras. Thus, $\psi^{-1}(\GG(x))$ is a Lie subalgebra
  which provides a direct sum decomposition of $\g_{2(2)}$ into irreducible $\sli(3,\R)$-modules. Such decomposition,
  by \cite[Chapter 22]{Fulton}, is given by $\psi^{-1}(\GG(x))\oplus\R^3\oplus\R^{3*}$. Therefore,
  \begin{equation*}
   \psi(\psi^{-1}(\GG(x)))\oplus\psi(\R^3)\oplus\psi(\R^{3*})=\GG(x)\oplus\psi(\R^3)\oplus\psi(\R^{3*}),
  \end{equation*}
  is a decomposition of $\HH$ into irreducible $\sli(3,\R)$-modules.

  Recall that we have a previous decomposition of $\HH$ into irreducible $\sli(3,\R)$-modules
  \begin{equation*}
   \HH=\GG(x)\oplus\VV(x)\oplus\VV^*(x).
  \end{equation*}
  Comparing the two decomposition of $\HH$ into irreducible $\sli(3,\R)$-modules we obtain our desired result.
 \end{proof}

 We fix an isomorphism of Lie algebras $\psi:\g_{2(2)}\to\HH$ as in the previous lemma. Let $G_{2(2)}$
 denote a simply connected Lie group such that Lie$(G_{2(2)})=\g_{2(2)}$. By Lemma \ref{lem: right act}, there exists
 an analytic isometric right $G_{2(2)}$-action on $\TM$ such that $\psi(X)=X^*$, for all $X\in\g_{2(2)}$. Now, we
 consider the map:
 \begin{eqnarray*}
  p:G_{2(2)}&\to&\TM\\
   g & \mapsto & x \cdot g,
 \end{eqnarray*}
 that satisfies $dp_e(X)=X^*_x=\psi(X)$ for every $X\in\g_{2(2)}$. Thus, by our choice of $\psi$ and Lemma
 \ref{lem: decomposition}, $dp_e$ is an isomorphism that maps $\sli(3,\R)$ on $T_x\FF$ and $\R^3\oplus\R^{3*}$
 onto $T_x\FF^\perp$. Since $p$ is $G_{2(2)}$-equivariant for the right action on its domain, then we have a
 local diffeomorphism.

 \begin{lemma}\label{lem: metric inv}
  Let $\bar{g}$ be the metric on $\g_{2(2)}$ defined as the pullback under $dp_e$ of the metric $g_x$ on
  $T_x\TM$ then $\bar{g}$ is $\sli(3,\R)$-invariant.
 \end{lemma}
 \begin{proof}

  By properties of $dp_e$ and the isomorphism $\psi$, we need only to prove that the metric on $\HH$ defined
  as the pullback of $g_x$ by the evaluation map, $X\mapsto X_x$, is $\GG(x)$-equivariant.

  Let $\tilde{g}$ be the metric on $\HH$ obtained of this way. Let $X,Y,Z\in\HH$ be given with $X\in\GG(x)$.
  By Lemma \ref{lem: evaluation} there
  exists $X_0\in\sli(3,\R)$ such that $X=\rho_x(X_0)+X_0^*$, where $\rho_x$ is the homomorphism in Proposition
  \ref{prop: g(x)} and $X_0^*$ is the vector field on $\TM$ induced by $X_0$ through the left $\TSL(3,\R)$-action.
  Therefore
  \begin{eqnarray*}
   \tilde{g}([X,Y],Z)&=&g_x([X,Y]_x,Z_x)=g([X,Y],Z)|_x=g([\rho_x(X_0)+X_0^*,Y],Z)|_x\\
    &=&g([\rho_x(X_0),Y],Z)|_x=\rho_x(X_0)g(Y,Z)|_x-g(Y,[\rho_x(X_0),Z])|_x\\
    &=&-g(Y,[\rho_x(X_0),Z])|_x=-g(Y,[\rho_x(X_0)+X_0^*,Z])|_x\\
    &=&-g(Y,[X,Z])|_x=-g_x(Y_x,[X_x,Z_x])=-\tilde{g}(Y,[X,Z]).
  \end{eqnarray*}
  Where we have used the fact that $\HH$ centralizes $X_0^*$ and that $\rho_x(X_0)$ is a Killing vector field,
  for the metric $g$, which vanishes in $x$. Thus, we take $\bar{g}$ as the pullback of $\tilde{g}$ by the
  isomorphism $\psi$ to obtain the desired result.
 \end{proof}

 Now, by the previous lemma and Lemma \ref{lem: uniqueproduct}, we can rescale the metric along the bundles $T\FF$
 and $T\FF^\perp$ in $\TM$ such that the new metric, $\widehat{g}$, on $\TM$ satisfies that
 $K=(dp_e)^*(\widehat{g}_x)$, is the Killing form on $\g_{2(2)}$.

 Since the elements of $\HH\subset\Kill(\TM)$ preserve the direct sum decomposition,
 $T\TM=T\FF\oplus T\FF^\perp$, then $\HH\subset\Kill(\TM,\widehat{g})$. Note that $\widehat{g}$ is invariant under
 both the left $\TSL(3,\R)$-action and the right $G_{2(2)}$-action on $\TM$. Observe, also, that the metric
 $\widehat{g}$ can be obtained from the lift of a correspondingly rescaled metric on $M$.

 Consider the bi-invariant metric on $G_{2(2)}$ induced by the Killing form $K$, which we denote with the same letter.
 The previous argument and discussion imply that the local diffeomorphism
 \begin{equation*}
  p:(G_{2(2)},K)\to(\TM,\widehat{g})
 \end{equation*}
 is a local isometry. This last property of $p$, the completeness of $(G_{2(2)},K)$ and the simply
 completeness of $(\TM,\widehat{g})$ imply, by Corollary 20 in \cite[p. 202]{Oneill}, that $p$ is an isometry.
 
 \begin{proposition}\label{prop: G_2(2)}
  Let $M$ be an analytic connected finite volume pseudo-Rie\-ma\-nnian manifold with $\dim(M)=14$. If $M$ is complete,
  admits an analytic and isometric right $\SL(3,\R)$-action with a dense orbit such that $\HH$
  (Lemma \ref{lem: embed}) is a simple Lie algebra with $\dim(\HH)=14$. Then there exists an analytic diffeomorphism
  $p:G_{2(2)}\to\TM$ and an analytic isometric right $G_{2(2)}$-action on $\TM$ such that:
  \begin{itemize}
   \item[(i)] On $\TM$, the left $\TSL(3,\R)$-action and the right $G_{2(2)}$-action commute,
   \item[(ii)] p is $G_{2(2)}$-equivariant for the right $G_{2(2)}$-action on its domain,
   \item[(iii)] for a pseudo-Riemannian metric $\widehat{g}$ in $\TM$ obtained by rescaling the original metric on
    the summands of the decomposition $T\TM=T\FF\oplus T\FF^\perp$, the map
    \begin{equation*}
     p:(G_{2(2)},K)\to(\TM,\widehat{g})
    \end{equation*}
    is an isometry where $K$ is the bi-invariant metric on $G_{2(2)}$ induced from the Killing form of its Lie
    algebra.
  \end{itemize}
 \end{proposition}

 Considering $G_{2(2)}$ with the bi-invariant pseudo-Riemannian metric $K$, induced by the Killing form of its Lie
 algebra, we can assume that $(G_{2(2)},K)$ is the isometric universal covering space of $(\TM,\widehat{g})$.

 By Proposition $4.5$ of \cite{OQ} we have that the isometry group of the pseudo-Rie\-man\-nian
 manifold $(G_{2(2)},K)$, which we denote by $\Iso(G_{2(2)})$, has only a finite number of connected components.
 Such proposition also shows that
 \begin{equation*}
  \Iso(G_{2(2)})_0=L(G_{2(2)})R(G_{2(2)}),
 \end{equation*}
 where $L(G_{2(2)})$ and $R(G_{2(2)})$ are the subgroups of left and right translations on $G_{2(2)}$, respectively.

 Let $\varrho:\TSL(3,\R)\to\text{Iso}(G_{2(2)})$ be the homomorphism induced by the isometric left
 $\TSL(3,\R)$-action on $G_{2(2)}$. From the previous observations, the covering
 \begin{equation*}
  G_{2(2)}\times G_{2(2)}\to L(G_{2(2)})R(G_{2(2)})
 \end{equation*}
 yields the existence of homomorphisms $\varrho_1,\varrho_2:\TSL(3,\R)\to G_{2(2)}$ such that
 \begin{equation*}
  \varrho(g)=L_{\varrho_1(g)}\circ R_{\varrho_2(g)^{-1}} \qquad \forall g\in\TSL(3,\R).
 \end{equation*}

 By Proposition \ref{prop: G_2(2)}, we have that $\varrho(g)\circ R_h=R_h\circ\varrho(g)$ for all
 $g\in\TSL(3,\R)$ and $h\in G_{2(2)}$, which implies that $\varrho_2(\TSL(3,\R))$ is
 contained in the center of $G_{2(2)}$, thence $\varrho(g)=L_{\varrho_1(g)}$ for all
 $g\in\TSL(3,\R)$. Thus, the $\TSL(3,\R)$-action on $G_{2(2)}$ is induced by the
 homomorphism $\varrho_1:\TSL(3,\R)\to G_{2(2)}$ and the left action of $ G_{2(2)}$ onto itself. Note,
 by our hypotheses, that the homomorphism $\varrho_1$ is non-trivial.

 By Proposition \ref{prop: G_2(2)} we have that $\pi_1(M)\subset\Iso(G_{2(2)})$ and by the previous
 observations $\Gamma_1=\pi_1(M)\cap\Iso_0( G_{2(2)})$ is a finite index subgroup of $\pi_1(M)$. Therefore,
 for every $\gamma\in\Gamma_1$ there exist $h_1,h_2\in G_{2(2)}$ such that $\gamma=L_{h_1}\circ R_{h_2}.$

 Since the left $\TSL(3,\R)$-action on $G_{2(2)}$ is the lift of an action on $M$,
 this left $\TSL(3,\R)$-action commutes with the $\Gamma_1$-action. Applying that    property to
 $L_{h_1}\circ R_{h_2}=\gamma\in\Gamma_1$
 we obtain that $L_{h_1}\circ L_{\varrho_1(g)}=L_{\varrho_1{(g)}}\circ L_{h_1}$ for all $g\in\TSL(3,\R)$, thus
 $\Gamma_1\in L(Z(\TSL(3,\R))) R( G_{2(2)})$ where $Z(\TSL(3,\R))$ is the centralizer
 of $\varrho_1(\TSL(3,\R))$ in $ G_{2(2)}$. By  Lemma \ref{lem: center}, the center of
 $G_{2(2)}$ has finite index in $Z(\TSL(3,\R))$ and therefore
 $R(G_{2(2)})$ has finite index in $L(\TSL(3,\R))R( G_{2(2)})$. In particular, $\Gamma=\Gamma_1\cap R( G_{2(2)})$
 is a finite index subgroup of $\Gamma_1$ and also of $\pi_1(M)$.

 The natural identification of $R(G_{2(2)})$ with $G_{2(2)}$ realizes $\Gamma$ as a discrete
 subgroup of $G_{2(2)}$ such that $G_{2(2)}/\Gamma$ is a finite covering space of $M$.

 Let $\varphi:G_{2(2)}/\Gamma\to M$ be the corresponding covering map. For the left
 $\widetilde{SL}(3,\R)$-action on $ G_{2(2)}/\Gamma$ given by the homomorphism
 $\varrho_1:\widetilde{SL}(3,\R)\to  G_{2(2)},$ the constructions in the previous paragraphs show that the map
 $\varphi$ is $\TSL(3,\R)$-equivariant. Finally, we note that $\varphi$ is a local isometry for the metric
 $\widehat{g}$, on $\TM$ considered in Proposition \ref{prop: G_2(2)}.

 Next, we are going to show that the subgroup $\Gamma$ is a lattice in $G_{2(2)}$. For the proof of that
 result it is enough to prove that $M$ has finite volume in the metric $\widehat{g}$. Recall, we are assuming
 that $M$ has finite volume in its original metric.

 \begin{lemma}\label{lem: volume}
  If \emph{vol} and $\emph{vol}_{\widehat{g}}$ denote the volume elements on $M$ for the original metric and the
  rescaled metric, respectively. Then, there is some constant $C>0$ such that
  $\emph{vol}_{\widehat{g}}=C\emph{vol}$.
 \end{lemma}
 \begin{proof}
  We consider $(x^1,x^2,\ldots,x^{14})$ some coordinate of $M$ in a neighborhood $U$ of a given point such that
  $(x^1,\ldots,x^8)$ defines a set of coordinates of the leaves of the foliation $\FF$ in such neighborhood.
  For the original metric $g$ on $M$, consider the orthogonal bundle $T\FF^\perp$ and a set of $1$-forms
  $\theta^1,\ldots,\theta^6$ that define a basis for its dual $(T\FF^\perp)^*$ at every point in $U$. Thus,
  in $U$ the metric $g$ has an expression of the form:
  \begin{equation*}
   g=\sum^8_{i,j=1}l_{ij}dx^i\otimes dx^j + \sum_{i,j=1}^6h_{ij}\theta^i\otimes\theta^j.
  \end{equation*}
  From this and the definition of the volume element as an $14$-form, we have:
  \begin{equation*}
   \text{vol}=\sqrt{|\det(l_{ij})\det(h_{ij})|}dx^1\wedge\ldots\wedge dx^8\wedge\theta^1\wedge\cdots\wedge\theta^6.
  \end{equation*}
  On the other hand, since the metric $\widehat{g}$ is obtained by rescaling $g$ along the bundles $T\OO$ and
  $T\OO^\perp$, then has an expression of the form:
  \begin{equation*}
   \widehat{g}=\sum^8_{i,j=1}C_1l_{ij}dx^i\otimes dx^j + \sum_{i,j=1}^6C_2h_{ij}\theta^i\otimes\theta^j.
  \end{equation*}
  for some constants $C_1,C_2\neq0$. Therefore, the volume element of $\widehat{g}$ satisfies:
  \begin{eqnarray*}
   \text{vol}_{\widehat{g}}&=&\sqrt{|\det(C_1l_{ij})\det(C_2h_{ij})|}dx^1\wedge\ldots\wedge dx^8
    \wedge\theta^1\wedge\cdots\wedge\theta^6 \\
    &=&\sqrt{|C_1^8C_2^6|}\text{vol}.
  \end{eqnarray*}
 \end{proof}

\subsubsection{$\HH\simeq\sli(4,\R)$}

 Here, we assume $\HH_0(x)\simeq\R$ and that $\HH$ is a simple Lie algebra then, by Lemma \ref{lem: H0=R}(2),
 our centralizer $\HH$ is isomorphic to $\sli(4,\R)$. Therefore, we have our following result
 \begin{lemma}\label{lem: isomorphism H=R}
  There is an isomorphism
  \begin{equation*}
   \psi:
   \sli(4,\R)=\sli(3,\R)\oplus\R\oplus\R^3\oplus\R^{3*}\to
   \GG(x)\oplus\HH_0(x)\oplus\VV(x)\oplus\VV^*(x)=\HH
  \end{equation*}
  of Lie algebras that preserves the summands in that order. In particular, $\psi$ is an isomorphism of
  $\sli(3,\R)$-modules.
 \end{lemma}
 \begin{proof}
   Similar to the proof of Lemma \ref{lem: iso to g22}.
 \end{proof}

 Let us fix an isomorphism of Lie algebras $\psi:\sli(4,\R)\to\HH$ as in the previous lemma. Let $\TSL(4,\R)$ denote
 the universal covering of $\SL(4,\R)$ then Lie$(\TSL(4,\R))=\sli(4,\R)$. By Lemma \ref{lem: right act}, there exists
 an analytic isometric right $\TSL(4,\R)$-action on $\TM$ such that $\psi(X)=X^*$, for every $X\in\sli(4,\R)$. This
 right action centralizes the left $\TSL(3,\R)$-action on $\TM$ and thus preserves the bundles $T\FF$ and
 $T\FF^\perp$.

 Given the previous right $\TSL(4,\R)$-action on $\TM$ we define the map $p:\TSL(4,\R)\to\TM$ defined as
 $g\mapsto x\cdot g$, for all $g\in\TSL(4,\R)$.
 Observe that this map is $\TSL(4,\R)$-invariant for the right action on its domain satisfying
 $dp_e(X)=X^*_x=\psi(X)$ for every $X\in\sli(3,\R)$. Note that $dp_e$ is surjective with
 $\ker(dp_e)=\psi^{-1}(\HH_0(x)$.

 Let $H$ be a connected subgroup of $\TSL(4,\R)$ such that $\text{Lie}(H)=\psi^{-1}(\HH_0(x))$. Here,
 $H$ is not a compact subgroup and, by exercise $(vi)$ in \cite[p. 152]{Helgason}, a closed subgroup. Hence,
 the map 
 \begin{equation*}
  \overline{p}:H\backslash\TSL(4,\R)\to\TM, \qquad [g]=Hg \mapsto x \cdot g,
 \end{equation*}
 for all $Hg\in H\backslash\TSL(4,\R)\to\TM$, is well defined. Observe that
 $T_{He}\big(H\backslash\TSL(4,\R)\big)=\sli(3,\R)\oplus\R^3\oplus\R^{3*}$.

 By our choice of the map $\psi$ we have that $d\overline{p}_{He}$ is an isomorphism which maps $\sli(3,\R)$ onto
 $T\FF$ and $\R^3\oplus\R^{3*}$ onto $T\FF^\perp$. Since $\overline{p}$ is an $\TSL(4,\R)$-equivariant map for
 the right action on its domain then $\overline{p}$ is an analytic local diffeomorphism at $He$.

 The definition of the map $\overline{p}$ and Lemma \ref{lem: isomorphism H=R} imply that
 that $d\overline{p}_{He}=ev_x\circ\psi$ restricted to the subspace $\sli(3,\R)\oplus\R^3\oplus\R^{3*}$.
 This map induces a metric on $H\backslash\TSL(4,\R)$, which is the result of the next lemma.

 \begin{lemma} \label{lem: pullback H=R}
  Let $\bar{g}$ be the metric on $T_{He}(H\backslash\TSL(4,\R))$ defined as the pullback under
  $d\bar{p}_{[e]}$ of the metric $g_x$ on $T_x\TM$, then $\bar{g}$ is $\sli(3,\R)$-invariant.
 \end{lemma}
 \begin{proof}
 Since $d\bar{p}_{[e]}=ev_x\circ\psi|_{\big(\sli(3,\R)\oplus\R^3\oplus\R^{3*}\big)}$, then $d\bar{p}_{He}$
 is a homomorphism
 of $\sli(3,\R)$-modules. Recall that the $\sli(3,\R)$-module structure in $\sli(3,\R)\oplus\R^3\oplus\R^{3*}$ is
 given by the subalgebra $\sli(3,\R)$ and in $T_x\TM$ by
 the subalgebra $\rho_x(\sli(3,\R)(x))$, as in Proposition \ref{prop: g(x)}(4).

 Since the metric $g$ in $T_x\TM$ is invariant under the action of $\rho_x(\sli(3,\R)(x))$ we have that
 if $u,v\in T_{He}(H\backslash\TSL(4,\R))$ and $X\in\sli(3,\R)$ then, by Lemma \ref{lem: evaluation}
 there are $U,V\in\HH$ such that such that $U_x=u$ and $V_x=v$ and by Proposition \ref{prop: g(x)}
 $\widehat{\rho}_x(X)\in\GG(x)$ and
 $X_0=\psi^{-1}(\widehat{\rho}_x(X))\in T_{He}(H\backslash\TSL(4,\R))$ satisfying that
 \begin{eqnarray*}
  \bar{g}([X,v],w)&:=&\bar{g}([X_0,v],w)=g_x(d\bar{p}_{He}([X_0,v]),d\bar{p}_{He}(w))\\
         &=&g_x((ev_x\circ\psi)[X_0,v],(ev_x\circ\psi)(w))=g_x(\psi([X_0,v])_x,\psi(w)_x)\\
         &=&g(\psi([X_0,v]),\psi(w))|_x=g([\psi(X_0),\psi(v)],\psi(w))|_x\\
         &=&g([\widehat{\rho}_x(X),V],W)|_x=g([\rho_x(X)+X^*,V],W)|_x\\
         &=&g([\rho_x(X),V],W)|_x=\rho_x(X)(g(V,W))|_x-g(V,[\rho_x(X),W])|_x\\
         &=&-g(V,[\rho_x(X),W])|_x=-g(V,[\rho_x(X)+X^*,W])|_x\\
         &=&-g(V,[\widehat{\rho}_x(X),W])|_x=-g(\psi(v),[\psi(X_0),\psi(w)])|_x\\
         &=&-g(\psi(v),\psi([X_0,w]))|_x=-g_x(d\bar{f}_{[e]}(v),d\bar{f}_{He}([X_0,w]))\\
         &=&-\bar{g}(v,[X_0,w])=-\bar{g}(v,[X,w]),
 \end{eqnarray*}
 Recall that $\HH$ centralizes $X^*$ and $\rho_{x_0}(X)$ is a Killing field for the metric $g$ in $\TM$.
 Therefore, by properties of the maps $ev_{x_0}$ and $\varphi$, the metric $\bar{g}$ in
 \begin{equation*}
  T_{He}(H\backslash\TSL(4,\R))=\sli(3,\R)\oplus\R^3\oplus\R^{3*}.
 \end{equation*} is $\sli(3,\R)$-invariant.
 \end{proof}

 Next, we analyze the pseudo-Riemannian metric structure of the analytic manifold $H\backslash\TSL(4,\R)$.

 First, let $K_n$ be the Killing form on $\sli(n,\R)$ with $n\geq2$. Recall that
 $K_n(X,Y)=2n\cdot\mathrm{tr}(XY)$ for all $X,Y\in\sli(n,\R)$.
 By the decomposition of $\sli(4,\R)$ as $\sli(3,\R)$-module in Section \ref{sect: Action} and
 the definition of the Killing form we have that $\sli(3,\R)$, $\R^3\oplus\R^{3*}$ and
 $\sli(3,\R)\oplus\R^3\oplus\R^{3*}$ (as subspaces of $\sli(4,\R)$) are non-degenerate with respect to the Killing
 form $K_4$.

 Denote by $K^1,K^2$ and $K$ the restriction of $K_4$ to the subspaces $\sli(3,\R)$, $\R^3\oplus\R^{3*}$ and
 $\sli(3,\R)\oplus\R^3\oplus\R^{3*}$, respectively. Since $K_4$ is invariant by the adjoint action of
 $\sli(4,\R)$ then it is clear that $K_4$ is also invariant under the adjoint action of $\sli(3,\R)$, by
 its inclusion in $\sli(4,\R)$. Therefore, $K^1$, $K^2$ and $K$ are invariant by the adjoint action
 of $\sli(3,\R)$.

 \begin{remark}\label{rem: multiple}
  By the definition and properties of $K^1$ and since $\sli(3,\R)$ is a simple Lie algebra we have, by Schur's Lemma,
  that $K^1=c_1K_3$, with $c_1\neq0$. On the other hand, by Section \ref{sect: Action} and Lemma
  \ref{lem: uniqueform}, we have that $K^2=K_4|_{\R^3\oplus\R^{3*}}$ is a non-zero multiple of a unique
  $\sli(3,\R)$-invariant symmetric bilinear form on $\R^3\oplus\R^{3*}$.
 \end{remark}

 Let $\pi:\TSL(4,\R)\to H\backslash\TSL(4,\R)$ be the natural quotient map. With respect to the previous map,
 we assume that
 \begin{equation}
  d\pi_0=d\pi|_{\sli(3,\R)\oplus\R^3\oplus\R^{3*}}:\sli(3,\R)\oplus\R^3\oplus\R^{3*}\to T_{He}H\backslash\TSL(4,\R) \label{eq: lin isometry}
 \end{equation} is a linear isometry.

 By its construction, we have that the quotient space $H\backslash\TSL(4,\R)$ is a \emph{reductive coset}
 (see \cite[p. 310]{Oneill}). On the other hand, by properties of the Killing form $K_4$ and the definition of $K$,
 we have that $K$ is Ad$(H)$-invariant. Thus, such result and the isometry
 \eqref{eq: lin isometry} imply, by Proposition 22 in \cite[p. 311]{Oneill}, the $\TSL(4,\R)$-invariance
 on $H\backslash\TSL(4,\R)$. Therefore, by Lemma \ref{lem: isomorphism H=R}, the properties of $K_4$ and the
 isometry \eqref{eq: lin isometry} we have that the manifold $H\backslash\TSL(4,\R)$ is a
 \emph{naturally reductive homogeneous space}, (see \cite[p. 312]{Oneill}).

 As a consequence that our manifold is a naturally reductive homogeneous space we have that $H\backslash\TSL(4,\R)$
 is complete (see \cite[p. 313]{Oneill}), hence, by Lemma 24 in \cite[p. 312]{Oneill}, our quotient map
 $\pi:\TSL(4,\R)\to H\backslash\TSL(4,\R)$ is a pseudo-Riemannian submersion.

 Next, we show that we can rescale the metric on $\TM$ such that the pullback of this new metric, with respect to
 the map $\bar{p}$, implies that \eqref{eq: lin isometry} is, effectively, a linear isometry. But,
 first we need the following result.

 \begin{lemma}\label{lem: rescale-metric H=R}
  Let $\LA\cdot,\cdot\RA_1$ and $\LA\cdot,\cdot\RA_2$ be the inner products on
  $\sli(3,\R)$ and ${\R^3}\oplus{\R^{3*}}$, respectively. Assume that
  $\LA\cdot,\cdot\RA_1$ and $\LA\cdot,\cdot\RA_2$ are $\sli(3,\R)$-invariant. Then there
  exist $c_1,c_2\in\R$ such that
  \begin{equation*}
   c_1\LA\cdot,\cdot\RA_1+c_2\LA\cdot,\cdot\RA_2,
  \end{equation*}
  is $K$, the Killing form of $\sli(4,\R)$ restricted to $\sli(3,\R)\oplus{\R^3}\oplus{\R^{3*}}$.
 \end{lemma}
 \begin{proof}
  Recall, Schur's Lemma implies that in $\g$, a simple real Lie algebra with a simple complexification, any
  $\g$-invariant non-degenerate symmetric bilinear form on $\g$ is a multiple by a real scalar of the Killing
  form.

  On the other hand, we have proved in Lemma \ref{lem: uniqueproduct} that there is, up to a multiple by a real scalar,
  a unique $\sli(3,\R)$-invariant non-degenerate bilinear form on $\R^3\oplus\R^{3*}$.

  Now, the result follows from previous results.
 \end{proof}

 \begin{remark}\label{rem: rescale}
  By the results in Lemmas \ref{lem: pullback H=R} and \ref{lem: rescale-metric H=R}, we can rescale the metric $g$
  along the bundles $T\FF$ and $T\FF^\perp$ on $\TM$ such that the new metric, $\widehat{g}$, satisfies
  $\big(d\bar{f}_{He}\big)^*(\widehat{g}_x)=K$, the Killing form on $\sli(4,\R)$ restricted to
  $\sli(3,\R)\oplus\R^3\oplus\R^{3*}$.
 \end{remark}

 Since the elements in $\HH$ preserve the decomposition $T\TM=T\FF\oplus T\FF^\perp$, then
 $\HH\subset\Kill(\TM,\widehat{g})$. That is, the elements in $\HH$ are Killing vector fields for the metric
 $\widehat{g}$ rescaled as in Remark \ref{rem: rescale}. Therefore, $\widehat{g}$ is invariant under the right
 $\TSL(4,\R)$-action. In a similar way, the left $\TSL(3,\R)$-action on $\TM$ preserves the rescaled metric
 $\widehat{g}$. Note that $\widehat{g}$ on $\TM$ can be obtained as the lift of a correspondingly rescaled
 metric $\widehat{g}$ in $M$.

 Remark \ref{rem: rescale} implies that the local diffeomorphism
 \begin{equation*}
  \bar{p}:(H\backslash\TSL(4,\R),K)\to(\TM,\widehat{g})
 \end{equation*} is a local isometry. Thus, the completeness of $(H\backslash\TSL(4,\R),K)$ and the simple
 completeness of $\TM$ imply, by Corollary 29 in \cite[p. 202]{Oneill}, that $\bar{p}$ is an isometry.
 Therefore our next result.

 \begin{proposition}\label{prop: almost H=R}
  Let $M$ be a connected analytic pseudo-Riemannian manifold. Suppose that $M$ is complete, has finite volume and
  admits an analytic and isometric $\SL(3,\R)$-action with a dense orbit such that the centralizer of this action,
  $\HH$, is a Lie algebra simple with $\dim(\HH)=15$. If $\dim(M)=14$, then there exists an analytic diffeomorphism
  $\bar{p}:H\backslash\TSL(4,\R)\to\TM$ and an analytic isometric right $\TSL(4,\R)$-action on $\TM$ such that:
  \begin{itemize}
   \item[(1)] on $\TM$ the left $\TSL(3,\R)$-action and the right $\TSL(4,\R)$-action commute with each other,
   \item[(2)] $\bar{p}$ is $\TSL(4,\R)$-equivariant for the right $\TSL(4,\R)$-action on its domain,
   \item[(3)] for a pseudo-Riemannian metric $\widehat{g}$ in $\TM$ obtained by rescaling the original metric
    on the summands of $T\TM=T\FF\oplus T\FF^\perp$, the map
    \begin{equation*}
    \bar{p}:(H\backslash\TSL(4,\R),K)\to(\TM,\widehat{g})
    \end{equation*}
    is an isometry where $K$ is the metric on $H\backslash\TSL(4,\R)$ which makes of the quotient map
    $\pi$ a pseudo-Riemannian submersion.
  \end{itemize}
 \end{proposition}

\appendix\section{}
 With the definitions and results from \cite{Vinberg} we have the next Lemma.

 \begin{lemma} \label{lem: max}
  Suppose that $\rho:\sli(3,\R)\to\g_{2(2)}$ is an injective Lie algebra homomorphism. Then $\s=\rho(\sli(3,\R))$ is a
  subalgebra of $\g_{2(2)}$ with its centralizer, $\mathfrak{z_g(s)}$, equal to zero.
 \end{lemma}

 \begin{lemma} \label{lem: center}
  Suppose that $G$ is a connected Lie group locally isomorphic to ${G}_{2(2)}$ and consider
  $\iota:\TSL(3,\R)\to G$ a non trivial homomorphism of Lie groups. Then, the centralizer
  $Z_G(\iota(\TSL(3,\R)))$ of $\iota(\TSL(3,\R))$ in $G$ contains $Z(G)$ (the center of $G$) as a finite index subgroup.
 \end{lemma}
 \begin{proof}
  Let $S=\iota(\TSL(3,\R))$ and denote the Lie algebra of $S$ as $\s$. Since $Z(G)\subseteq Z_G(S)$ and,
  by the previous Lemma, $\mathfrak{z_g(s)}=0$ then $Z_G(S)$ and $Z(G)$ are discrete. The proof that $Z_G(S)$
  is finite is a consequence of Lemma 1.1.3.7 in \cite{Warner}.
 \end{proof}

  We prove now the relationship between the Killing form on the simple Lie group $\g_{2(2)}$ and the
  $\sli(3,\R)$-invariant bilinear form, both obtained in $\sli(3,\R)$ and $\R^3\oplus\R^{3*}$. But before we
  need the following result.

 \begin{lemma}\label{lem: uniqueform}
  There is, up to a multiple by a real scalar, exactly one $\sli(3,\R)$-invariant non-degenerate bilinear form on
  $\R^3\oplus\R^{3*}$.
 \end{lemma}
 \begin{proof}
  We have proved that there exists in $\R^3\oplus\R^{3*}$ a $\sli(3,\R)$-invariant non-degenerate bilinear form.
  This has as consequence the existence of an isomorphism $\varrho:\R^3\oplus\R^{3*}\to(\R^3\oplus\R^{3*})^*$ of
  $\sli(3,\R)$-modules.

  Then we have an isomorphism $\varrho(\C):\C^3\oplus\C^{3*}\to\C^3\oplus\C^{3*}$ of
  $\sli(3,\C)$-modules, that by Schur's Lemma, is just
  the multiple of the identity by a complex number when res\-tricted to $\C^3$ and to another complex number when
  res\-tricted to $\C^{3*}$. Furthermore, since $\varrho(\C)$ is the complexification of $\varrho$ we have that
  these numbers are real.

  The result follows from the previous arguments and the fact that $\R^3$ and $\R^{3*}$ belong to the nullcone
  of the inner product.
 \end{proof}

 \begin{lemma}\label{lem: uniqueproduct}
  Let $\LA\cdot,\cdot\RA_1$ and $\LA\cdot,\cdot\RA_2$ be inner products on $\sli(3,\R)$ and
  $\R^3\oplus\R^{3*}$, respectively. If we suppose that $\LA\cdot,\cdot\RA_1$ and
  $\LA\cdot,\cdot\RA_2$ are $\sli(3,\R)$-invariant. Then, there exist $a_1,a_2\in\R$ such that
  $a_1\LA\cdot,\cdot\RA_1+a_2\LA\cdot,\cdot\RA_2$ is the Killing form of $\g_{2(2)}$.
 \end{lemma}
 \begin{proof}
  Recall that Schur's Lemma implies that in $\g$, a simple real Lie algebra with a simple complexification, any
  $\g$-invariant non-degenerate symmetric bilinear form on $\g$ is unique up to a multiple, that is, the multiple by
  a scalar of the Killing form.

  In particular, we have that $\LA\cdot,\cdot\RA_1$ is a multiple of the Killing form of $\g_{2(2)}$ ($K$) when
  restricted to $\sli(3,\R)$, this is $\LA X,Y\RA_1=c_1K|_{\sli(3,\R)}(X,Y)$ for all $X,Y\in\sli(3,\R)$ and
  some non-zero $c_1\in\R$.

  On the other hand, from Lemma \ref{lem: uniqueform} we have the existence of a non-zero scalar $c_2\in\R$ such that
  $\LA\cdot,\cdot\RA_2=c_2K|_{\R^3\oplus\R^{3*}}(\cdot,\cdot)$.

  Now, the result is a consequence of the previous arguments and Lemma \ref{lem: uniqueform}.
 \end{proof}

\end{document}